\documentclass{article}
\usepackage{geometry}               
\usepackage{tikz}
\usetikzlibrary{shapes.geometric, arrows}

\usepackage{graphicx}
\usepackage{caption}
\usepackage{subcaption}
\usepackage{comment}
\usepackage{amsmath} % assumes amsmath package installed
\usepackage{amssymb}  % assumes amsmath package installed
\usepackage{amsthm}  % assumes amsmath package installed
\usepackage{bm}
\usepackage{lipsum}
\usepackage{color}
\usepackage{multirow}
\usepackage{array}
\usepackage{cite}
\usepackage{hyperref}
\usepackage{comment}

\usepackage{algorithm}
\usepackage{algpseudocode}
\usepackage{enumitem}
\usepackage{float}
\usepackage{url}

\newtheorem{prop}{Proposition}

\newtheorem{lem}{Lemma}
\newtheorem{thm}{Theorem}
\newtheorem{rem}{Remark}
\newtheorem{assumption}{Assumption}

\numberwithin{equation}{section}

\newcommand{\R}{\mathbb{R}}

\newcommand{\de}{\mathrm{d}}

\DeclareMathOperator*{\argmax}{arg max}
\DeclareMathOperator*{\argmin}{arg min}

\title{Solving nonconvex Hamilton--Jacobi--Isaacs equations with PINN-based policy iteration}
\author{
    Hee Jun Yang\thanks{These authors contributed equally to this work.}\ \thanks{National Institute for Mathematical Sciences (NIMS), 70, Yuseong-daero 1689 beon-gil, Yuseong-gu, Daejeon 34047, Republic of Korea.}
    \and
    Minjung Gim\footnotemark[1]\ \footnotemark[2]
    \and
    Yeoneung Kim\thanks{Corresponding author. Department of Industrial Engineering, Yonsei University, 50 Yonsei-ro, Seodaemun-gu, Seoul 03722, Republic of Korea. E-mail: \texttt{yeoneung@yonsei.ac.kr}.}
}
\begin{document}

\maketitle

\begin{abstract}
We propose a mesh-free policy iteration framework that combines classical dynamic programming with physics-informed neural networks (PINNs) to solve high-dimensional, nonconvex Hamilton--Jacobi--Isaacs (HJI) equations arising in stochastic differential games and robust control. The method alternates between solving linear second-order PDEs under fixed feedback policies and updating the controls via pointwise minimax optimization using automatic differentiation. Under standard Lipschitz and uniform ellipticity assumptions, we establish local uniform convergence of the value function iterates to the unique viscosity solution, together with a rigorous global $L^2$ convergence rate. Numerical experiments demonstrate the accuracy and scalability of the method. In a two-dimensional stochastic path-planning game with a moving obstacle, our method matches finite-difference benchmarks with relative $L^2$-errors below $10^{-2}$. In high-dimensional publisher-subscriber differential games with anisotropic noise, the proposed approach consistently outperforms direct PINN solvers. Our results suggest that integrating PINNs with policy iteration is a practical and theoretically grounded method for solving high-dimensional, nonconvex HJI equations, with potential applications in robotics, finance, and multi-agent reinforcement learning.
\end{abstract}


\section{Introduction}
The Hamilton--Jacobi--Isaacs (HJI) equation plays a fundamental role in differential games and robust control, characterizing the value function of zero-sum stochastic dynamic games. In the presence of diffusion, the HJI equation takes the form
\[
\partial_t v + H(t,x,\nabla_x v) = -\tfrac12 \operatorname{Tr}(\sigma\sigma^\top D_{xx}^2 v),
\]
with suitable terminal condition. The Hamiltonian $H$ typically involves a minimax operation over control variables and may be nonconvex or nonsmooth, which presents substantial analytical and numerical challenges.

Classical numerical methods, such as finite difference and semi-Lagrangian schemes~\cite{kushner1990numerical,souganidis1985approximation,barles1991convergence}, provide robust convergence guarantees via viscosity solution theory, but their reliance on structured spatial meshes renders them intractable in high dimensions due to the curse of dimensionality. Various extensions using unstructured meshes~\cite{barth1998numerical}, convexification~\cite{Klibanov2021}, and radial basis collocation~\cite{cecil2004numerical,Ferretti2016} have been proposed, but scalability remains a limiting factor.

To alleviate grid-related bottlenecks, mesh-free methods based on physics-informed neural networks (PINNs)~\cite{raissi2019physics,han2018solving} have emerged as promising alternatives. PINNs approximate solutions of PDEs by minimizing residuals through neural networks with automatic differentiation, and have been applied to Hamilton--Jacobi (HJ) and Hamilton--Jacobi--Bellman (HJB) equations. Recent studies~\cite{Liu2022} demonstrate improved convergence in nonconvex HJ settings using adaptive losses. However, directly minimizing nonconvex residuals can be unstable and prone to poor local minima, especially when the Hamiltonian involves a nonsmooth minimax structure.

To address these limitations, policy iteration (PI) methods have been introduced for HJB/HJI equations~\cite{kerimkulov2020exponential,guo2025policy,kawecki2022discontinuous,tran2025policy}. By alternating between value evaluation and policy improvement steps, PI offers a structured approach that improves stability and convergence. When coupled with PINNs, this framework enables mesh-free iteration and can mitigate difficulties associated with nonconvexity in the Hamiltonian. In particular, Lee and Kim~\cite{lee2025hamilton} proposed a PI framework based on deep operator learning (DeepONet) to solve Hamilton--Jacobi--Bellman equations, proposing the first operator-learning-based implementation of policy iteration for HJB equations. Their approach emphasized function-space policy updates via neural operators and showed promising results in high-dimensional optimal control problems. Our method departs from this line by adopting a residual-based PINN formulation tailored to nonconvex HJI equations, while providing a rigorous convergence guarantee under ellipticity assumptions.

In this paper, we propose a PINN-based policy iteration framework for solving nonconvex HJI equations. At each step, the value function is approximated by a neural network trained to minimize the PDE residual for fixed policies. Gradients from automatic differentiation are used to update feedback controls via pointwise optimization. This leads to a continuous optimization framework that avoids spatial grids entirely. We establish convergence to the viscosity solution under standard Lipschitz assumptions. In addition to empirical performance, the proposed method offers a structural advantage: it allows for explicit control over approximation errors via provable $L^2$ bounds, something rarely feasible in nonconvex HJI problems.

The rest of the paper is organized as follows. In Section~\ref{sec:prelim}, we briefly review the formulation of Hamilton--Jacobi--Isaacs equations in the context of differential games. Section~\ref{sec:method} presents the proposed PINN-based policy iteration framework, detailing the policy evaluation and improvement steps, the associated training procedure, and a theoretical analysis of the method, including a convergence guarantee under suitable regularity and ellipticity conditions. Section~\ref{sec:exp} presents the results of extensive numerical experiments on benchmark problems in multiple dimensions, highlighting the accuracy, scalability, and robustness of the proposed approach. Finally, Section~\ref{sec:conclusion} concludes the paper and discusses potential directions for future research.

\section{Stochastic Differential Games and Hamilton--Jacobi--Isaacs equations}\label{sec:prelim}
Given $0\leq t<T$, $d,m_1,m_2 \in\mathbb{N}$, $A\subset \mathbb{R}^{m_1}$ and $B\subset \mathbb{R}^{m_2}$, we consider a two-player zero-sum stochastic differential game over a finite horizon $ [0, T] $, where the state process $ X(s) \in \mathbb{R}^d $ evolves according to the controlled stochastic differential equation
\begin{equation} \label{eq:sde}
	\begin{cases}
        \begin{aligned}
    		\mathrm{d} X(s) &= f(s,X(s), a(s), b(s)) \mathrm{d}s + \sigma(s,X(s))\mathrm{d}W_s, \quad \text{for}\quad s \geq t,\\
    		X(t) &= x,
        \end{aligned}
	\end{cases}
\end{equation}
where $ W_s \in \mathbb{R}^d $ denotes a standard $ d $-dimensional Brownian motion and the functions $ f : [0,T]\times \mathbb{R}^d \times A \times B \to \mathbb{R}^d $ and $ \sigma : [0,T]\times \mathbb{R}^d \to \mathbb{R}^{d \times d} $ represent the drift and diffusion coefficients, respectively. Let $\mathcal{F}:=(\mathcal{F}_s)_{s \geq 0}$ be the filtration generated by $(W_s)_{s\geq 0}$, and for $s\in[t,T]$, set $a \in \mathcal{A}_t$ and $b \in \mathcal{B}_t$ where
\[
\begin{cases}
    \begin{aligned}
    	\mathcal{A}_t &=  \{ a : [t,T] \to A \mid a \text{ is a } \{\mathcal{F}_s\}_{s \in [t,T]} \text{-adapted process}  \}, \\
    	\mathcal{B}_t &=  \{ b : [t,T] \to B \mid b \text{ is a } \{\mathcal{F}_s\}_{s \in [t,T]} \text{-adapted process}  \}.
    \end{aligned}
\end{cases}
\]
The performance of a control pair $ (a, b) \in \mathcal{A}_t \times \mathcal{B}_t $ is evaluated through the cost functional
\begin{equation*}
	J(t, x; a, b) = \mathbb{E} \left [ \int_t^T c(s,X(s), a(s), b(s)) \mathrm{d}s + g(X_T)\right  ],
\end{equation*}
where $ c : \mathbb{R}^d \times A \times B \to \mathbb{R} $ is the running cost, and $ g : \mathbb{R}^d \to \mathbb{R} $ is the terminal cost. Player I seeks to minimize this expected cost, while Player II aims to maximize it.  The set of admissible nonanticipative strategies for Player II beginning at time $t$ is defined by
\[
\Gamma_t = \{ \beta : \mathcal{A}_t \to \mathcal{B}_t \mid \text{nonanticipating} \},
\]
where the strategy $\beta$ is called nonanticipative if, for \( a_1, a_2 \in \mathcal{A}_t \) and \( s \in [t, T] \),
\[
a_1(\cdot) = a_2(\cdot) \text{ on } [t, s) \quad \Longrightarrow \quad \beta[a_1](\cdot) = \beta[a_2](\cdot) \text{ on } [t, s).
\]

Accordingly, the value function of the game is defined as
\begin{equation*}
	v(t, x) = \sup_{\beta \in \Gamma_t} \inf_{a \in \mathcal{A}_t}J(t,x; a, \beta[a]).
\end{equation*}

Under standard regularity assumptions, such as Lipschitz continuity and boundedness of $ f $, $c$, and $\sigma$ in both $t$ and $x$, and uniform ellipticity of the matrix $\sigma \sigma^\top(t,x) $, the value function $ v(t,x) $ satisfies the dynamic programming principle and is characterized as the unique viscosity solution of the Hamilton--Jacobi--Isaacs (HJI) equation
\[
\begin{cases}
	\partial_t v(t,x) + H(t,x,\nabla_x v(t,x)) = -\frac{1}{2} \operatorname{Tr}(\sigma\sigma^\top(t,x) D_{xx}^2 v(t,x)) \quad &\textnormal{in}\quad [0,T) \times \mathbb{R}^d,\\
	v(T,x)=g(x) \quad&\text{on}\quad \mathbb{R}^d,
\end{cases}
\]
where 
\[
H(t, x, p) := \sup_{\textbf{b} \in B} \inf_{\textbf{a} \in A} L(t, x, p)(\textbf{a}, \textbf{b})
\]
with
\[
L(t, x, p)(\textbf{a}, \textbf{b}) := c(t, x, \textbf{a}, \textbf{b}) + p \cdot f(t, x, \textbf{a}, \textbf{b}).
\]

Each term in the HJI equation admits a natural interpretation. The term $\partial_t v$ describes the temporal evolution of the value, which propagates backward from the terminal data $v(T,\cdot)=g$. The Hamiltonian $H(t,x,\nabla_x v)$ encodes the instantaneous minimax interaction between the two players: the Lagrangian $L$ aggregates the running cost $c$ and the directional change of the value along the controlled drift, $\nabla_x v \cdot f$, and the $\sup$/$\inf$ structure selects each player's optimal instantaneous response given the marginal value $\nabla_x v$. Finally, the second-order term $\frac{1}{2}\operatorname{Tr}(\sigma\sigma^\top D_{xx}^2 v)$, which arises from It\^o's formula applied to the state dynamics \eqref{eq:sde}, accounts for the diffusive effect of the noise and regularizes the equation under the uniform ellipticity assumption.

\section{Physics-informed approach for solving HJI equations}\label{sec:method}
\subsection{Policy iteration for HJI equations}
We begin by introducing the notation used throughout the paper. For \( x \in \mathbb{R}^d \), we write \( |x| \) for the Euclidean norm. Given a function \( f: \Omega \to \mathbb{R}^n \), we denote its standard \( L^p \) norm by
\[
\|f\|_p := \left( \int_\Omega |f|^p  \mathrm{d}x\right)^{1/p} \quad \text{for}\quad p \in (0, \infty],
\]
where \( |f| \) denotes the pointwise Euclidean norm of \( f \). We say \( f \in C^{k,\beta}(\Omega) \) for \( k \in \mathbb{N} \) and \( \beta \in (0,1) \) if all partial derivatives \( D^\alpha f \) of order \( |\alpha| \leq k \) exist and are continuous on \( \Omega \), and for all multi-indices \( \alpha \) with \( |\alpha| = k \),
\[
[D^\alpha f]_{C^{0,\beta}(\Omega)} := \sup_{x \neq y \in \Omega} \frac{|D^\alpha f(x) - D^\alpha f(y)|}{|x - y|^\beta} < \infty.
\]
We write \( C^\beta(\Omega) := C^{0,\beta}(\Omega) \) for the space of \(\beta\)-Hölder continuous functions, and \( C(\Omega) \) for the space of continuous functions.

We now introduce the policy iteration framework used throughout the paper. Our goal is to solve high-dimensional Hamilton--Jacobi--Isaacs (HJI) equations by combining policy iteration with physics-informed neural networks (PINNs). The proposed method alternates between solving linear PDEs under fixed feedback policies and updating the feedback controls via gradient-based minimax steps.

As a starting point, we follow the iterative scheme introduced in~\cite{guo2025policy}, where a discrete-time policy iteration (PI) algorithm for the HJI equation is first proposed based on the mesh-free algorithm, which is demonstrated in Algorithm~\ref{alg:PIstep}.
\begin{algorithm}[H]
	\caption{Mesh-free Policy Iteration for HJI}\label{alg:PIstep}
	\begin{algorithmic}[1]
		\Statex \textbf{Input:} Lipschitz continuous initial feedback $(\alpha_0,\beta_0)$
		\For{$n = 0,1,2,\dots$}
		\State \textbf{Policy evaluation:} find $v_n$ solving
		\[
		\partial_t v_n + L(t,x,\nabla_x v_n)(\alpha_n,\beta_n)
		= -\tfrac12\mathrm{Tr}(\sigma\sigma^\top D_{xx}^2 v_n), \quad v_n(T,x) = g(x).
		\]
		\State \textbf{Policy improvement:} at each $(t,x)$ with gradient $p = \nabla_x v_n(t,x)$,
		\[
		\begin{aligned}
			\alpha_{n+1,b}(t,x) &\gets \argmin_{\mathbf{a}\in A} L(t,x,p)(\mathbf{a},b), \\
			\beta_{n+1}(t,x) &\gets \argmax_{\mathbf{b}\in B} L(t,x,p)(\alpha_{n+1,b}(t,x), \mathbf{b}), \\
			\alpha_{n+1}(t,x) &\gets \alpha_{n+1,\beta_{n+1}(t,x)}(t,x).
		\end{aligned}
		\]
		\EndFor
		\Statex \textbf{Output:} converged feedback pair $(\alpha_n,\beta_n)$ and value $v_n$
	\end{algorithmic}
\end{algorithm}
We emphasize that the policy improvement step in Algorithm~\ref{alg:PIstep} is a nested optimization computing the saddle point of $(\mathbf a,\mathbf b)\mapsto L(t,x,p)(\mathbf a,\mathbf b)$: for each fixed $\mathbf b\in B$, $\alpha_{n+1,b}$ is Player~I's best response; $\beta_{n+1}$ is Player~II's optimal strategy against these best responses; and the final policy $\alpha_{n+1}$ is Player~I's best response to $\beta_{n+1}$. Under Assumption~\ref{assump:main}(i) below, this sequential procedure is well-defined and returns the unique saddle point of $L(t,x,p)$; see Lemma~\ref{lem:lip}.

In their setting, the diffusion matrix $\sigma$ may be degenerate, and the value function $v$ can fail to be differentiable, making the feedback-based policy update
\[
(\textbf{a},\textbf{b})\mapsto\argmin_{\textbf{a}\in A}\argmax_{\textbf{b}\in B} L(t,x,\nabla_x v)(\textbf{a},\textbf{b})
\]
ill-posed. To address this, the authors of~\cite{guo2025policy} introduced a discrete space-time grid and an artificial viscosity term that ensures enough regularity to define the minimax update at grid points. Convergence to the viscosity solution is then obtained via a careful limiting argument.

In contrast, uniform ellipticity provides sufficient parabolic regularity for the policy-evaluation iterates so that their spatial gradients admit continuous representatives. Together with the single-valued continuous minimax selector, this makes the feedback update well defined pointwise. We emphasize that the convergence analysis below uses only uniform boundedness of the induced drifts, rather than Lipschitz continuity of the feedback policies.

We recall the form of the Lagrangian
\[
L(t, x, p)(\textbf{a}, \textbf{b}) := c(t, x, \textbf{a}, \textbf{b}) + p \cdot f(t, x, \textbf{a}, \textbf{b}), 
\]
and introduce assumptions on the dynamics and cost function.
\begin{assumption}\label{assump:main}
	We impose the following assumptions throughout the paper.
	\begin{enumerate}[label=(\roman*), leftmargin=*]
		\item The control sets $A \subset \mathbb{R}^{m_1}$ and $B \subset \mathbb{R}^{m_2}$ are convex and compact. For each fixed $(t,x,p)\in [0,T]\times \mathbb{R}^d \times \mathbb{R}^d$, the map
		\[
		(\mathbf{a},\mathbf{b}) \mapsto L(t,x,p)(\mathbf{a},\mathbf{b}) := c(t,x,\mathbf{a},\mathbf{b}) + p\cdot f(t,x,\mathbf{a},\mathbf{b})
		\]
		is $\mu_A$-strongly convex in $\mathbf{a}$ and $\mu_B$-strongly concave in $\mathbf{b}$.
		
		That is, for every $(t,x,p)$, the maps
		\[
		\mathbf{a} \mapsto L(t,x,p)(\mathbf{a}, \mathbf{b}) \quad\text{and}\quad
		\mathbf{b} \mapsto L(t,x,p)(\mathbf{a}, \mathbf{b})
		\]
		are strongly convex and strongly concave, respectively, uniformly in $(t,x,p)$.
		
		\item The functions $f$, $c$, $g$, and $\sigma$ are bounded and Lipschitz continuous in all variables, with common Lipschitz constant $L_u>0$. In addition, there exists $\beta > 0$ such that:
		\begin{enumerate}[label=(\alph*), leftmargin=*]
			\item The terminal cost satisfies $g\in C_b^{2+\beta}(\mathbb R^d)\cap L^2(\mathbb R^d).$;
			\item The running cost satisfies
			\[
            \|\,\sup_{(\mathbf a,\mathbf b)\in A\times B}|c(\cdot,\cdot,\mathbf a,\mathbf b)|\,\|_{L^\infty(0,T;L^2(\mathbb{R}^d))} < \infty
			\]
			\item The mappings $(t,x) \mapsto f(t,x,\mathbf{a},\mathbf{b})$ and $\sigma(t,x)$ belong to $C^\beta([0,T]\times\mathbb{R}^d)$, uniformly in $(\mathbf{a},\mathbf{b})$.
		\end{enumerate}
		\item The diffusion coefficient $ \sigma \sigma^\top \in C([0,T] \times \mathbb{R}^d; \mathbb{R}^{d \times d}) $ is uniformly elliptic: there exists $ \lambda > 0 $ such that
		\[
		\sigma\sigma(t,x)^\top  \succeq \lambda I_d \quad \text{for all}\quad (t,x) \in [0,T] \times \mathbb{R}^d.
		\]

	\end{enumerate}
\end{assumption}
The assumptions stated above lead to a mild structural property of the policy-update map; this property, formulated below, is central to our convergence analysis.

\begin{lem}[Lipschitz continuity of the feedback selector]\label{lem:lip}
	Let the Lagrangian $L(t, x, p)(\textbf{a}, \textbf{b}) := c(t, x, \textbf{a}, \textbf{b}) + p \cdot f(t, x, \textbf{a}, \textbf{b})$. Then for every $(t,x,p)$, the policy update
	\begin{equation*}
		\begin{cases}
            \begin{aligned}
			\alpha^\star(t,x,p) &=
			\argmin_{\mathbf a\in A}L(t,x,p)(\mathbf a,\beta^\star(t,x,p)),\\
			\beta^\star(t,x,p) &=
			\argmax_{\mathbf b\in B}L(t,x,p)(\alpha^\star(t,x,p),\mathbf b),
            \end{aligned}
		\end{cases}
	\end{equation*}
	admits a unique solution, and $(\alpha^\star,\beta^\star)$ is globally Lipschitz in~$p$:
	\[
	|\alpha^\star(t,x,p_1)-\alpha^\star(t,x,p_2)|
	+
	|\beta^\star(t,x,p_1)-\beta^\star(t,x,p_2)|
	\le
	\kappa |p_1-p_2|,
	\]
	with a constant~$\kappa:=\frac{2L_u}{\mu}>0$.
\end{lem}

\begin{proof}
Fix $(t,x)$ and write
\[
L_i(a,b):=L(t,x,p_i)(a,b),
\qquad
z_i:=(a_i,b_i):=z^\star(t,x,p_i),
\quad i=1,2.
\]
Since $z_i$ is the saddle point of $L_i$, the strong convexity in $a$
and strong concavity in $b$ imply
\[
L_1(a_2,b_1)-L_1(a_1,b_2)
\ge
\frac{\mu_A}{2}|a_1-a_2|^2
+
\frac{\mu_B}{2}|b_1-b_2|^2
\]
and
\[
L_2(a_1,b_2)-L_2(a_2,b_1)
\ge
\frac{\mu_A}{2}|a_1-a_2|^2
+
\frac{\mu_B}{2}|b_1-b_2|^2.
\]
Adding these inequalities and recalling that
\[
L_i(a,b)=c(t,x,a,b)+p_i\cdot f(t,x,a,b),
\]
we obtain
\[
\begin{aligned}
\mu_A|a_1-a_2|^2+\mu_B|b_1-b_2|^2
&\le
(p_1-p_2)\cdot
\bigl(
f(t,x,a_2,b_1)-f(t,x,a_1,b_2)
\bigr)\\
&\le
L_u|p_1-p_2|
\bigl(|a_1-a_2|+|b_1-b_2|\bigr),
\end{aligned}
\]
where $L_u$ is the Lipschitz constant of $f$ in the control variables.
Setting $\mu:=\min\{\mu_A,\mu_B\}$ and using
\[
|a_1-a_2|^2+|b_1-b_2|^2
\ge
\frac12
\bigl(|a_1-a_2|+|b_1-b_2|\bigr)^2,
\]
we conclude that
\[
|a_1-a_2|+|b_1-b_2|
\le
\frac{2L_u}{\mu}|p_1-p_2|.
\]
\end{proof}

We first record a compactness property of the exact policy-evaluation
iterates. This will be combined below with a weighted energy estimate
to identify the limit and obtain a quantitative convergence rate.

\begin{lem}[Uniform bounds and local compactness of the policy iterates]
\label{lem:PI_compactness}
Suppose Assumption~\ref{assump:main} holds, and let
$\{v_n\}_{n\ge0}$ be generated by Algorithm~\ref{alg:PIstep}.
Then there exists a constant $M>0$, independent of $n$, such that
\[
\sup_{n\ge0}
\|v_n\|_{L^\infty([0,T]\times\mathbb R^d)}
\le M.
\]
Moreover, for every $R>0$, there exist
$\alpha\in(0,1)$ and $C_R>0$, independent of $n$, such that
\[
\|v_n\|_{C^{\alpha/2,\alpha}([0,T]\times B_R)}
\le C_R.
\]
Consequently, $\{v_n\}_{n\ge0}$ is relatively compact in
$C_{\mathrm{loc}}([0,T]\times\mathbb R^d)$.
\end{lem}

\begin{proof}
Set
\[
b_n(t,x)
:=
f\bigl(t,x,\alpha_n(t,x),\beta_n(t,x)\bigr),
\qquad
c_n(t,x)
:=
c\bigl(t,x,\alpha_n(t,x),\beta_n(t,x)\bigr).
\]
Since the control sets are compact and the coefficients $f$ and $c$
are uniformly bounded,
\[
\|b_n\|_\infty\le \|f\|_\infty,
\qquad
\|c_n\|_\infty\le \|c\|_\infty
\]
uniformly in $n$.

For each $n$, the Feynman--Kac representation gives
\[
v_n(t,x)
=
\mathbb E\left[
g(X_T^{t,x})
+
\int_t^T c_n(s,X_s^{t,x})\,\de s
\right],
\]
where
\[
\de X_s^{t,x}
=
b_n(s,X_s^{t,x})\,\de s
+
\sigma(s,X_s^{t,x})\,\de W_s,
\qquad
X_t^{t,x}=x.
\]
It follows that
\[
|v_n(t,x)|
\le
\|g\|_\infty
+
(T-t)\|c\|_\infty,
\]
and hence
\[
\sup_{n\ge0}\|v_n\|_\infty
\le
\|g\|_\infty+T\|c\|_\infty
=:M.
\]

Each $v_n$ solves the uniformly parabolic linear equation
\[
\partial_t v_n
+
\frac12\operatorname{Tr}
\bigl(\sigma\sigma^\top D_{xx}^2v_n\bigr)
+
b_n\cdot\nabla_xv_n
=
-c_n,
\qquad
v_n(T,\cdot)=g.
\]
The ellipticity constants, the $L^\infty$ bounds of $b_n$ and $c_n$,
and the coefficient regularity of $\sigma\sigma^\top$ are all uniform
in $n$. Therefore, standard interior and terminal-boundary Hölder
estimates for uniformly parabolic equations imply that, for every
$R>0$, there exist $\alpha\in(0,1)$ and $C_R>0$, independent of $n$,
such that
\[
\|v_n\|_{C^{\alpha/2,\alpha}([0,T]\times B_R)}
\le C_R.
\]
The Arzel\`a--Ascoli theorem and a diagonal argument now imply the
relative compactness of $\{v_n\}_{n\ge0}$ in
$C_{\mathrm{loc}}([0,T]\times\mathbb R^d)$.
\end{proof}

The following theorem is the uniformly elliptic continuous-space
analogue of \cite[Theorem~1.1]{guo2025policy}. In addition to identifying
the limit as the HJI solution, it gives a global exponential rate in
$C([0,T];L^2(\mathbb R^d))$.

\begin{thm}[Convergence and uniform exponential rate]
\label{thm:PI_uniform}
\label{prop:exp_rate}
Fix $T>0$ and suppose Assumption~\ref{assump:main} holds.
Let $\{v_n\}_{n\ge0}$ be generated by Algorithm~\ref{alg:PIstep}.
Then $\{v_n\}$ converges locally uniformly on
$[0,T]\times\mathbb R^d$ to the unique bounded continuous viscosity
solution $v$ of
\begin{equation}\label{eq:hji}
\begin{cases}
\partial_t v
+
H(t,x,\nabla_xv)
=
-\dfrac12\operatorname{Tr}
\bigl(\sigma\sigma^\top(t,x)D_{xx}^2v\bigr),
&
\textnormal{in }[0,T)\times\mathbb R^d,
\\[0.4em]
v(T,x)=g(x),
&
\textnormal{on }\mathbb R^d.
\end{cases}
\end{equation}
Moreover, there exist constants $C>0$ and $\rho\in(0,1)$ such that
\[
\sup_{t\in[0,T]}
\|v_n(t,\cdot)-v(t,\cdot)\|_{L^2(\mathbb R^d)}
\le
C\rho^n
\]
for every $n\ge0$.
\end{thm}

\begin{proof}
Throughout the proof, we write
\[
\delta_n:=v_{n+1}-v_n,
\qquad
\pi_n:=(\alpha_n,\beta_n),
\qquad
p_n:=\nabla_xv_n.
\]
Recall from Algorithm~\ref{alg:PIstep} that, for each $n\ge0$, the iterate $v_n$ is the solution of the linear policy-evaluation equation
\begin{equation}\label{eq:policy_eval}
\partial_t v_n + L(t,x,\nabla_x v_n)(\alpha_n,\beta_n)
= -\tfrac12\operatorname{Tr}\bigl(\sigma\sigma^\top D_{xx}^2 v_n\bigr),
\qquad v_n(T,x)=g(x).
\end{equation}
After reversing time, we work on $[0,T]\times\mathbb R^d$ with
initial data prescribed at $t=0$, so that
\[
\delta_n(0,\cdot)=0
\qquad
\text{for every }n\ge0.
\]

\smallskip
\emph{Step 1: Uniform $L^2$ bounds and contraction of successive
differences.}

Define
\[
c_n(t,x)
:=
c\bigl(t,x,\alpha_n(t,x),\beta_n(t,x)\bigr),
\]
and
\[
\bar c(t,x)
:=
\sup_{(\mathbf a,\mathbf b)\in A\times B}
|c(t,x,\mathbf a,\mathbf b)|.
\]
Then
\[
\|c_n(t,\cdot)\|_2
\le
\|\bar c(t,\cdot)\|_2
\]
uniformly in $n$, and
\[
\bar c\in
L^\infty(0,T;L^2(\mathbb R^d))
\]
by Assumption~\ref{assump:main}(ii)(b).

Applying Proposition~\ref{prop:L2} to each policy-evaluation equation \eqref{eq:policy_eval}
gives
\begin{equation}\label{eq:M0_bound}
\begin{split}
\sup_{n\ge0}
\left[
\sup_{t\in[0,T]}\|v_n(t,\cdot)\|_2^2
+
\lambda
\int_0^T\|\nabla_xv_n(t,\cdot)\|_2^2\,\de t
\right]
\\
\le
C_T
\left(
\|g\|_2^2
+
T
\|\bar c\|_{L^\infty(0,T;L^2)}^2
\right)
=:M_0<\infty.
\end{split}
\end{equation}
In particular,
\[
v_n\in
C([0,T];L^2(\mathbb R^d))
\cap
L^2(0,T;H^1(\mathbb R^d))
\]
uniformly in $n$.

Subtracting the policy-evaluation equations \eqref{eq:policy_eval} for $v_{n+1}$ and $v_n$,
we obtain
\begin{equation}\label{eq:delta_PDE}
\partial_t\delta_n
-
\frac12
\operatorname{Tr}
\bigl(\sigma\sigma^\top D_{xx}^2\delta_n\bigr)
=
\mathcal L_n,
\qquad
\delta_n(0,\cdot)=0,
\end{equation}
where
\[
\begin{split}
\mathcal L_n
&:=
L(t,x,p_{n+1})(\pi_{n+1})
-
L(t,x,p_n)(\pi_n)
\\
&=
\underbrace{
L(t,x,p_{n+1})(\pi_{n+1})
-
L(t,x,p_n)(\pi_{n+1})
}_{=:\mathrm I}
\\
&\quad+
\underbrace{
L(t,x,p_n)(\pi_{n+1})
-
L(t,x,p_n)(\pi_n)
}_{=:\mathrm{II}}.
\end{split}
\]
Since $L$ is affine in $p$ with coefficient $f$,
\[
|\mathrm I|
\le
\|f\|_\infty|\nabla_x\delta_n|.
\]

By the policy-improvement step,
\[
H(t,x,p_n)
=
L(t,x,p_n)(\pi_{n+1}),
\]
and, for $n\ge1$,
\[
H(t,x,p_{n-1})
=
L(t,x,p_{n-1})(\pi_n).
\]
Consequently,
\[
\begin{aligned}
|\mathrm{II}|
&=
\left|
H(t,x,p_n)-H(t,x,p_{n-1})
+
L(t,x,p_{n-1})(\pi_n)
-
L(t,x,p_n)(\pi_n)
\right|
\\
&\le
2\|f\|_\infty|p_n-p_{n-1}|.
\end{aligned}
\]
Thus, for every $n\ge1$,
\begin{equation}\label{eq:Rn_pointwise}
|\mathcal L_n|
\le
\|f\|_\infty
\left(
|\nabla_x\delta_n|
+
2|\nabla_x\delta_{n-1}|
\right).
\end{equation}

Let
\[
A(t,x):=\sigma\sigma^\top(t,x).
\]
Testing \eqref{eq:delta_PDE} with $\delta_n$ and using
\[
\operatorname{Tr}(AD^2w)
=
\operatorname{div}(A\nabla w)
-
(\operatorname{div}A)\cdot\nabla w,
\]
we obtain
\begin{equation}\label{eq:energy}
\begin{split}
\frac{\de}{\de t}
\|\delta_n(t,\cdot)\|_2^2
+
\lambda
\|\nabla_x\delta_n(t,\cdot)\|_2^2
&\le
2C_\sigma
\int_{\mathbb R^d}
|\nabla_x\delta_n||\delta_n|\,\de x
\\
&\quad+
2
\int_{\mathbb R^d}
|\mathcal L_n||\delta_n|\,\de x,
\end{split}
\end{equation}
where
\[
C_\sigma
:=
\frac12
\|\operatorname{div}(\sigma\sigma^\top)\|_\infty.
\]
Here $\lambda>0$ is the uniform ellipticity constant of Assumption~\ref{assump:main}(iii); it enters through the coercivity bound $\int_{\mathbb R^d}(\sigma\sigma^\top\nabla_x\delta_n)\cdot\nabla_x\delta_n\,\de x\ge\lambda\|\nabla_x\delta_n\|_2^2$ and quantifies the parabolic dissipation that drives the contraction below.
Substituting \eqref{eq:Rn_pointwise} and applying Young's inequality,
we find that, for every $\varepsilon>0$,
\begin{equation}\label{eq:diff_ineq}
\frac{\de}{\de t}
\|\delta_n(t,\cdot)\|_2^2
+
\frac{3\lambda}{4}
\|\nabla_x\delta_n(t,\cdot)\|_2^2
\le
K
\|\delta_n(t,\cdot)\|_2^2
+
\varepsilon
\|\nabla_x\delta_{n-1}(t,\cdot)\|_2^2,
\end{equation}
where
\[
K
:=
\frac{4(C_\sigma+\|f\|_\infty)^2}{\lambda}
+
\frac{4\|f\|_\infty^2}{\varepsilon}.
\]

Introduce the weighted quantities
\[
\widetilde E_n
:=
\sup_{t\in[0,T]}
e^{-Kt}
\|\delta_n(t,\cdot)\|_2^2,
\]
and
\[
\widetilde F_n
:=
\int_0^T
e^{-Kt}
\|\nabla_x\delta_n(t,\cdot)\|_2^2\,\de t.
\]
Multiplying \eqref{eq:diff_ineq} by $e^{-Kt}$ gives
\[
\frac{\de}{\de t}
\left(
e^{-Kt}
\|\delta_n(t,\cdot)\|_2^2
\right)
+
\frac{3\lambda}{4}
e^{-Kt}
\|\nabla_x\delta_n(t,\cdot)\|_2^2
\le
\varepsilon
e^{-Kt}
\|\nabla_x\delta_{n-1}(t,\cdot)\|_2^2.
\]
Integrating over $[0,t]$ and using $\delta_n(0,\cdot)=0$ yields
\[
e^{-Kt}
\|\delta_n(t,\cdot)\|_2^2
+
\frac{3\lambda}{4}
\int_0^t
e^{-Ks}
\|\nabla_x\delta_n(s,\cdot)\|_2^2\,\de s
\le
\varepsilon\widetilde F_{n-1}.
\]
Hence,
\begin{equation}\label{eq:contraction}
\widetilde E_n
\le
\varepsilon\widetilde F_{n-1},
\qquad
\widetilde F_n
\le
\frac{4\varepsilon}{3\lambda}
\widetilde F_{n-1}.
\end{equation}

Choose
\[
\varepsilon:=\frac{3\lambda}{8}.
\]
Then
\[
\widetilde F_n
\le
\frac12\widetilde F_{n-1},
\]
and therefore
\[
\widetilde F_n
\le
2^{-n}\widetilde F_0.
\]
Moreover,
\[
\widetilde E_n
\le
\frac{3\lambda}{4}
2^{-n}\widetilde F_0.
\]
By \eqref{eq:M0_bound},
\[
\widetilde F_0
\le
2
\int_0^T
\left(
\|\nabla_xv_1\|_2^2
+
\|\nabla_xv_0\|_2^2
\right)\,\de t
\le
\frac{4M_0}{\lambda}.
\]
Removing the exponential weight, we conclude that
\begin{equation}\label{eq:delta_decay}
\sup_{t\in[0,T]}
\|\delta_n(t,\cdot)\|_2
\le
C_0\rho^n,
\qquad
\rho:=\frac1{\sqrt2},
\end{equation}
where
\[
C_0
:=
\left(
\frac{3\lambda}{4}
e^{KT}\widetilde F_0
\right)^{1/2}.
\]

It follows that, for $m>n\ge1$,
\[
\begin{aligned}
\sup_{t\in[0,T]}
\|v_m(t,\cdot)-v_n(t,\cdot)\|_2
&\le
\sum_{k=n}^{m-1}
\sup_{t\in[0,T]}
\|\delta_k(t,\cdot)\|_2
\\
&\le
\frac{C_0}{1-\rho}\rho^n.
\end{aligned}
\]
Thus, $\{v_n\}$ is Cauchy in
$C([0,T];L^2(\mathbb R^d))$ and converges to some
\[
\bar v\in C([0,T];L^2(\mathbb R^d)).
\]
Furthermore,
\[
\sum_{n=0}^\infty
\|\nabla_x\delta_n\|_{L^2(0,T;L^2)}
<\infty,
\]
so that
\begin{equation}\label{eq:H1_strong_convergence}
v_n
\longrightarrow
\bar v
\quad\text{strongly in}\quad
C([0,T];L^2(\mathbb R^d))
\cap
L^2(0,T;H^1(\mathbb R^d)).
\end{equation}

\smallskip
\emph{Step 2: Identification of the limit.}

Define the policy-consistency defect
\[
r_n(t,x)
:=
L(t,x,\nabla_xv_n)(\pi_n)
-
H(t,x,\nabla_xv_n).
\]
Since $\pi_n$ is obtained by improving $v_{n-1}$,
\[
L(t,x,\nabla_xv_{n-1})(\pi_n)
=
H(t,x,\nabla_xv_{n-1}).
\]
Using the Lipschitz continuity of $L$ and $H$ in the gradient variable,
we obtain
\[
\begin{aligned}
|r_n|
&\le
\left|
L(t,x,\nabla_xv_n)(\pi_n)
-
L(t,x,\nabla_xv_{n-1})(\pi_n)
\right|
\\
&\quad+
\left|
H(t,x,\nabla_xv_{n-1})
-
H(t,x,\nabla_xv_n)
\right|
\\
&\le
2\|f\|_\infty
|\nabla_xv_n-\nabla_xv_{n-1}|
\\
&=
2\|f\|_\infty
|\nabla_x\delta_{n-1}|.
\end{aligned}
\]
Consequently,
\begin{equation}\label{eq:defect_convergence}
r_n
\longrightarrow0
\quad\text{strongly in}\quad
L^2(0,T;L^2(\mathbb R^d)).
\end{equation}

Returning to the original time orientation, the policy-evaluation
equation can be written as
\[
\partial_t v_n
+
\frac12\operatorname{Tr}
\bigl(\sigma\sigma^\top D_{xx}^2v_n\bigr)
+
H(t,x,\nabla_xv_n)
=
-r_n,
\qquad
v_n(T,\cdot)=g.
\]
Since
\[
A=\sigma\sigma^\top\in W^{1,\infty},
\]
we may write
\[
\operatorname{Tr}(AD^2w)
=
\operatorname{div}(A\nabla w)
-
(\operatorname{div}A)\cdot\nabla w
\]
and pass to the limit in the variational formulation.
By \eqref{eq:H1_strong_convergence}, the Lipschitz continuity of $H$
in $p$, and \eqref{eq:defect_convergence},
\[
H(t,x,\nabla_xv_n)
\longrightarrow
H(t,x,\nabla_x\bar v)
\quad\text{strongly in}\quad
L^2(0,T;L^2(\mathbb R^d)).
\]
Therefore, $\bar v$ is a variational solution of
\[
\partial_t\bar v
+
\frac12\operatorname{Tr}
\bigl(\sigma\sigma^\top D_{xx}^2\bar v\bigr)
+
H(t,x,\nabla_x\bar v)
=
0,
\qquad
\bar v(T,\cdot)=g.
\]

By Lemma~\ref{lem:PI_compactness}, there exists a subsequence
$\{v_{n_k}\}$ converging locally uniformly to a bounded continuous
function $w$. Since the same subsequence converges to $\bar v$ in
$C([0,T];L^2(\mathbb R^d))$, we have $w=\bar v$ almost everywhere on
every compact cylinder. We henceforth identify $\bar v$ with this
bounded continuous representative.

As shown above, $\bar v$ is a variational solution of
\[
\partial_t\bar v
+\frac12\operatorname{Tr}
\bigl(\sigma\sigma^\top D_{xx}^2\bar v\bigr)
+H(t,x,\nabla_x\bar v)=0,
\qquad
\bar v(T,\cdot)=g.
\]
By the standard weak--viscosity compatibility for uniformly parabolic
equations with Lipschitz gradient dependence
(cf.~\cite{ishii1995equivalence}), $\bar v$ is a viscosity solution of
\eqref{eq:hji}. The comparison principle therefore yields
\[
\bar v=v.
\]

\smallskip
\emph{Step 3: Local uniform convergence and conclusion.}

Let $\{v_{n_j}\}$ be an arbitrary subsequence. By
Lemma~\ref{lem:PI_compactness}, it admits a further subsequence
converging locally uniformly to some bounded continuous function
$\widetilde v$. Since the full sequence converges to $v$ in
$C([0,T];L^2(\mathbb R^d))$, we necessarily have $\widetilde v=v$.
Therefore,
\[
v_n\longrightarrow v
\quad\text{locally uniformly on}\quad
[0,T]\times\mathbb R^d.
\]

Finally, letting $m\to\infty$ in the Cauchy estimate gives
\[
\sup_{t\in[0,T]}
\|v_n(t,\cdot)-v(t,\cdot)\|_2
\le
\frac{C_0}{1-\rho}\rho^n,
\qquad n\ge1.
\]
The case $n=0$ follows by enlarging the constant $C$.
This completes the proof.
\end{proof}

Having established the convergence and exponential rate of the idealized policy iteration, we now turn to its practical PINN implementation.

\subsection{PINN-Based Policy Iteration}\label{subsec:pi}
In the practical implementation, each value-function iterate is
represented by a neural network $v_n(t,x;\theta_n)$. Its spatial
gradients and Hessians are computed by automatic differentiation,
so the pointwise policy-improvement step is well-defined throughout
the training domain.

For the fixed pair $ (\alpha_{n},\beta_{n}) $ the value function
is represented by a neural network
$ v_{n}(t,x;\theta_{n}) $.
Let
\(
\{(t_{j},x_{j})\}_{j=1}^{N_{\text{int}}}\subset(0,T)\times\mathbb{R}^{d}
\)
be interior collocation points and
\(
\{x^{T}_{k}\}_{k=1}^{N_{\text{bc}}}\subset\mathbb{R}^{d}
\)
terminal points.
The network parameters are obtained by minimizing

\begin{equation}\label{eq:loss}
	\begin{split}
		\mathcal{J}(\theta_{n})=
		&\frac1{N_{\text{int}}}\sum_{j=1}^{N_{\text{int}}}
		|
		\partial_{t}v_{n}+L(t_{j},x_{j},\nabla_x v_{n})(\alpha_{n},\beta_{n})
		+\tfrac12\operatorname{Tr}  (\sigma \sigma^\top(t_j,x_j)  D_{xx}^{2}v_{n}  )
		|^{2}
		\\\quad&+\frac1{N_{\text{bc}}}\sum_{k=1}^{N_{\text{bc}}}
		|v_{n}(T,x^{T}_{k})-g(x^{T}_{k})|^{2}.
	\end{split}
\end{equation}

When implemented, we parameterize the value function $v_n(t, x; \theta_n)$ using the following ansatz:
\begin{equation}\label{eq:ansatz}
	v_n(t, x; \theta_n) = g(x) + (T - t) \mathcal{N}_n(t, x; \theta_n),
\end{equation}
which explicitly enforces the terminal condition as a hard constraint, eliminating the need for a separate terminal loss term and thereby improving training stability and convergence speed~\cite{lagaris1998artificial}.

Since $\mathcal N_n$ is a finite feedforward network with sine
activations, for every fixed parameter vector $\theta_n$ the network
and its derivatives of any fixed order are globally bounded. In
particular,
\[
\mathcal N_n,\qquad
\partial_t\mathcal N_n,\qquad
\nabla_x\mathcal N_n,\qquad
D_{xx}^2\mathcal N_n
\in
L^\infty([0,T]\times\mathbb R^d).
\]
Consequently, under Assumption~\ref{assump:main}, the neural ansatz
$\tilde v_n$ is a bounded classical function with bounded derivatives
up to the order entering the PDE. Its residual $\mathcal R_n$ is
therefore bounded and continuous. The additional condition
\[
\mathcal R_n\in L^2(0,T;L^2(\mathbb R^d)),
\]
imposed in Theorem~\ref{thm:global_error}, controls its behavior at
spatial infinity.

All differential operators are evaluated by automatic differentiation.
Because no spatial grid is required, the procedure scales to high
state dimensions without suffering from the curse of dimensionality. 

\begin{algorithm}[H]
	\caption{Mesh-free PINN Policy Iteration (PINN-PI)}\label{algo:main}
	\begin{algorithmic}[1]
		\Require collocation sets $\{(t_j, x_j)\}$, terminal set $\{x_k^T\}$, tolerance ($\mathrm{tol})$, number of policy updates $M$, number of training epochs per policy update $E$
		\State choose any Lipschitz continuous initial feedback pair $(\alpha_{0},\beta_{0}) : [0,T]\times\mathbb{R}^{d}\to A\times B$
		\State initialize network parameters $\theta_0$ using Xavier initialization; set all biases to zero
		\For{$n = 0,1,2,\dots,M-1$}
		\For{$\ell = 1,\dots,E$} \Comment{\textbf{Value evaluation}}
		\State update $\theta_n$ of $v_n(t,x;\theta_n)$ via gradient descent to minimize $\mathcal{J}(\theta_n)$
		\If{$\ell \equiv 0 \pmod{100}$}
		\State resample $\{(t_j, x_j)\}$ uniformly from the computational domain.
		\EndIf
		\EndFor
		\ForAll{collocation points $(t,x)$ used in the gradient computation} \Comment{\textbf{Policy update}}
		\State compute $\nabla_x v_{n}(t,x;\theta_n)$ via automatic differentiation
		\State $\displaystyle\alpha_{n+1,b}(t,x)\gets\argmin_{\textbf{a}\in A} L  (t,x,\nabla_x v_{n}(t,x;\theta_n))(\textbf{a},b)$
		\State $\displaystyle\beta_{n+1}(t,x)\gets\argmax_{\textbf{b}\in B} L(t,x,\nabla_x v_{n}(t,x))(\alpha_{n+1,b}(t,x),\textbf{b})$
		\State $\alpha_{n+1}(t,x)\gets \alpha_{n+1,\beta_{n+1}(t,x)}(t,x)$
		\EndFor
        \If{$n>0$ \textbf{and} $\max_{(t,x) \in \mathcal{S}_{\text{test}}} | v_{n}(t,x) - v_{n-1}(t,x) | < \mathrm{tol}$} \Comment{\textbf{Optional}}
        \State \textbf{break}
        \EndIf
		\State next iteration by setting $\theta_{n+1} \gets \theta_n$ \Comment{\textbf{Warm-start}}
		\EndFor
	\end{algorithmic}
\end{algorithm}

%%%%%%%%%%%%%%%%%%%%%%%%%%%%%%%%%%%%%%%%%%%%%%%%%%%%%%%%%%%%%%%%%%%%%%%%
% Additional assumptions governing the policy-update consistency
%%%%%%%%%%%%%%%%%%%%%%%%%%%%%%%%%%%%%%%%%%%%%%%%%%%%%%%%%%%%%%%%%%%%%%%%

The implementation of the above policy iteration scheme relies on two key ingredients: (i) the ability to approximate value functions using neural networks trained on residual losses, and (ii) the ability to compute policy updates based on gradients obtained via automatic differentiation. Algorithm~\ref{algo:main} outlines the resulting PINN-based scheme, where each policy update is executed in a completely mesh-free setting, leveraging the continuous feedback gradient $\nabla_x v_n(t,x)$. 

When the collocation batch is refreshed during inner PINN training, the frozen
feedback laws are evaluated on the newly sampled points using the stored
value-network parameters that generated the current policies. This evaluation
does not constitute a policy-improvement step, but only recomputes the values
of the fixed feedback laws on the refreshed collocation batch.

Early stopping in Algorithm~\ref{algo:main} is optional; in our experiments we run a fixed number of outer iterations $M$ for consistency across runs. In practice, the outer iteration may be terminated once the difference between successive value iterates becomes small, as measured by a norm of $v_n-v_{n-1}$ on a fixed held-out set $\mathcal S_{\mathrm{test}}$ (e.g., $\ell^\infty$ or $\ell^2$).
Here, $\mathcal S_{\mathrm{test}}$ consists of time-state samples drawn from the computational domain and not used for training.

\begin{rem}[Comparison of Algorithm~\ref{alg:PIstep} and Algorithm~\ref{algo:main}]
	Algorithm~\ref{alg:PIstep} describes an idealized policy iteration framework under which each value function \( v_n \) is assumed to solve a linear elliptic--parabolic PDE exactly, using classical analytic or grid-based numerical methods. This setting enables a rigorous convergence analysis under viscosity solution theory but is limited to low-dimensional problems due to the need for exact PDE solvers. In contrast, Algorithm~\ref{algo:main} implements a practical mesh-free version based on physics-informed neural networks (PINNs). The value function \( v_n \) is approximated via a neural network trained by minimizing the PDE residual over sampled collocation points, and the gradient \( \nabla_x v_n \) required for policy updates is obtained via automatic differentiation. This structure enables high-dimensional scalability and smooth feedback policy updates, at the expense of introducing approximation error. %Such errors are rigorously controlled under Assumption~\ref{assump:policy_lip} and quantified through the residual tolerance and convergence analysis in Theorem~\ref{thm:global_error}.
\end{rem}

Since the practical algorithm employs neural network approximations and finite iterations, a theoretical justification of numerical errors is necessary. The next theorem quantifies how the total error can be controlled in terms of the residual at each step.

%%%%%%%%%%%%%%%%%%%%%%%%%%%%%%%%%%%%%%%%%%%%%%%%%%%%%%%%%%%%%%%%%%%%%%%%
% Main estimate incorporating policy-evaluation \emph{and} policy-update errors
%%%%%%%%%%%%%%%%%%%%%%%%%%%%%%%%%%%%%%%%%%%%%%%%%%%%%%%%%%%%%%%%%%%%%%%%
%%%%%%%%%%%%%%%%%%%%%%%%%%%%%%%%%%%%%%%%%%%%%%%%%%%%%%%%%%%%%%%%%%%%%%%%
% Main estimate incorporating policy-evaluation \emph{and} policy-update errors
%%%%%%%%%%%%%%%%%%%%%%%%%%%%%%%%%%%%%%%%%%%%%%%%%%%%%%%%%%%%%%%%%%%%%%%%

\begin{thm}[Global error of the practical PINN-PI algorithm]
\label{thm:global_error}
Suppose Assumption~\ref{assump:main} holds, and let
$\{\tilde v_n,\tilde\alpha_n,\tilde\beta_n\}_{n\ge0}$ be generated by
Algorithm~\ref{algo:main} with the sine-network
ansatz~\eqref{eq:ansatz}, where each $\mathcal N_n$ is a finite
feedforward network with sine activations.
Assume that the practical and exact schemes are initialized at the same
feedback pair,
\[
(\tilde\alpha_0,\tilde\beta_0)=(\alpha_0,\beta_0).
\]
Let $\mathcal R_n$ denote the corresponding PDE residual,
\[
\mathcal R_n
:=
\partial_t\tilde v_n
+
L(t,x,\nabla_x\tilde v_n)(\tilde\alpha_n,\tilde\beta_n)
+
\frac12
\operatorname{Tr}
\bigl(\sigma\sigma^\top D_{xx}^2\tilde v_n\bigr).
\]
Assume that
\[
\mathcal R_n
\in
C_b([0,T]\times\mathbb R^d)
\cap
L^2(0,T;L^2(\mathbb R^d))
\qquad
\text{for every }n\ge0,
\]
and set
\[
p_n
:=
\max_{0\le k\le n}
\|\mathcal R_k\|_{L^2(0,T;L^2(\mathbb R^d))}.
\]
Then there exist constants $C>0$ and $\rho\in(0,1)$, depending only on
the structural constants and the data norms appearing in
Assumption~\ref{assump:main}, such that
\begin{equation}\label{eq:global-error-bias}
\sup_{t\in[0,T]}
\|\tilde v_n(t,\cdot)-v(t,\cdot)\|_2
\le
C\bigl(p_n+\rho^n\bigr),
\end{equation}
where $v$ is the unique viscosity solution of the HJI
equation~\eqref{eq:hji}.
\end{thm}

\begin{proof}
Let $\{v_n,\pi_n\}_{n\ge0}$ denote the exact policy-iteration sequence
of Algorithm~\ref{alg:PIstep}, where
\[
\pi_n:=(\alpha_n,\beta_n),
\qquad
\tilde\pi_n:=(\tilde\alpha_n,\tilde\beta_n).
\]
For each $n\ge0$, let $\hat v_n$ be the exact solution of the linear
policy-evaluation equation obtained by freezing the practical feedback
law $\tilde\pi_n$. We decompose
\[
\tilde v_n-v
=
\underbrace{\tilde v_n-\hat v_n}_{=:A_n}
+
\underbrace{\hat v_n-v_n}_{=:B_n}
+
\underbrace{v_n-v}_{=:C_n}.
\]
We also write
\[
\Phi_n:=\tilde v_n-v_n=A_n+B_n,
\qquad
\delta_n:=v_{n+1}-v_n.
\]

By Lemma~\ref{lem:lip}, both schemes use the same single-valued feedback
selector. In particular, for every $n\ge1$,
\begin{equation}\label{eq:selector_consistency}
\tilde\pi_n
=
(\alpha^\star,\beta^\star)
(t,x,\nabla_x\tilde v_{n-1}),
\qquad
\pi_n
=
(\alpha^\star,\beta^\star)
(t,x,\nabla_xv_{n-1}).
\end{equation}
The common initialization gives $\hat v_0=v_0$.

\smallskip
\emph{Step 1: policy-evaluation error.}

Set
\[
b_n(t,x):=f(t,x,\tilde\pi_n(t,x)).
\]
Let $w_n$ be the unique variational solution of
\begin{equation}\label{eq:variational_residual_problem}
\begin{cases}
\displaystyle
\partial_t w_n
+
\frac12
\operatorname{Tr}
\bigl(\sigma\sigma^\top D_{xx}^2w_n\bigr)
+
b_n\cdot\nabla_xw_n
=
\mathcal R_n,
&
\text{in }[0,T)\times\mathbb R^d,
\\[0.4em]
w_n(T,\cdot)=0,
&
\text{on }\mathbb R^d.
\end{cases}
\end{equation}
Standard variational theory gives
\[
w_n
\in
C([0,T];L^2(\mathbb R^d))
\cap
L^2(0,T;H^1(\mathbb R^d)),
\qquad
\partial_tw_n
\in
L^2(0,T;H^{-1}(\mathbb R^d)).
\]
Since $\|b_n\|_\infty\le\|f\|_\infty$, Proposition~\ref{prop:L2}
yields
\begin{equation}\label{eq:wn_energy}
\sup_{t\in[0,T]}
\|w_n(t,\cdot)\|_2^2
+
\lambda
\int_0^T
\|\nabla_xw_n(t,\cdot)\|_2^2\,\de t
\le
C_A
\|\mathcal R_n\|_{L^2(0,T;L^2)}^2
\le
C_Ap_n^2.
\end{equation}

We next identify $w_n$ with the actual neural policy-evaluation error
\[
A_n:=\tilde v_n-\hat v_n.
\]
Indeed, subtracting the frozen-policy equation for $\hat v_n$ from the
residual equation satisfied by $\tilde v_n$ shows that $A_n$ solves
\begin{equation}\label{eq:An_residual_problem}
\begin{cases}
\displaystyle
\partial_t A_n
+
\frac12
\operatorname{Tr}
\bigl(\sigma\sigma^\top D_{xx}^2A_n\bigr)
+
b_n\cdot\nabla_xA_n
=
\mathcal R_n,
&
\text{in }[0,T)\times\mathbb R^d,
\\[0.4em]
A_n(T,\cdot)=0,
&
\text{on }\mathbb R^d.
\end{cases}
\end{equation}
The sine-network ansatz implies that $\tilde v_n$ is bounded and
continuous. Moreover, the frozen-policy solution $\hat v_n$ is bounded
and continuous by the Feynman--Kac representation. Hence $A_n$ is a
bounded continuous viscosity solution of
\eqref{eq:An_residual_problem}.

On the other hand, since $\mathcal R_n$ is bounded and continuous, the
variational solution $w_n$ admits a bounded continuous representative
and coincides with the viscosity solution of
\eqref{eq:variational_residual_problem}. The comparison principle for
this linear uniformly parabolic equation therefore gives
\[
A_n=w_n.
\]
Consequently,
\[
A_n
\in
C([0,T];L^2(\mathbb R^d))
\cap
L^2(0,T;H^1(\mathbb R^d)),
\qquad
\partial_tA_n
\in
L^2(0,T;H^{-1}(\mathbb R^d)),
\]
and \eqref{eq:wn_energy} becomes
\begin{equation}\label{eq:An_energy}
\sup_{t\in[0,T]}
\|A_n(t,\cdot)\|_2^2
+
\lambda
\int_0^T
\|\nabla_xA_n(t,\cdot)\|_2^2\,\de t
\le
C_Ap_n^2.
\end{equation}
\smallskip
\emph{Step 2: policy-mismatch equation.}

Subtracting the policy-evaluation equations for $\hat v_n$ and $v_n$
gives
\begin{equation}\label{eq:Bn_PDE_explicit}
\begin{split}
\partial_tB_n
&+
\frac12
\operatorname{Tr}
\bigl(\sigma\sigma^\top D_{xx}^2B_n\bigr)
+
f(t,x,\tilde\pi_n)\cdot\nabla_xB_n
\\
&=
-\mathcal Q_n,
\qquad
B_n(T,\cdot)=0,
\end{split}
\end{equation}
where
\[
\mathcal Q_n
:=
L(t,x,\nabla_xv_n)(\tilde\pi_n)
-
L(t,x,\nabla_xv_n)(\pi_n).
\]

We estimate $\mathcal Q_n$ without using a uniform
$L^\infty$-bound for $\nabla_xv_n$. For readability, fix $(t,x)$ and
set
\[
p:=\nabla_xv_n(t,x),
\qquad
q:=\nabla_x\tilde v_{n-1}(t,x),
\qquad
r:=\nabla_xv_{n-1}(t,x).
\]
By the policy-improvement identities,
\[
L(t,x,q)(\tilde\pi_n)=H(t,x,q),
\qquad
L(t,x,r)(\pi_n)=H(t,x,r).
\]
Consequently,
\[
\begin{aligned}
\mathcal Q_n
={}&
L(t,x,p)(\tilde\pi_n)
-
L(t,x,q)(\tilde\pi_n)
\\
&+
H(t,x,q)-H(t,x,r)
\\
&+
L(t,x,r)(\pi_n)
-
L(t,x,p)(\pi_n).
\end{aligned}
\]
Since both $L(t,x,\cdot)(\pi)$ and $H(t,x,\cdot)$ are globally
Lipschitz in the gradient variable with Lipschitz constant
$\|f\|_\infty$, we obtain
\[
\begin{aligned}
|\mathcal Q_n|
&\le
\|f\|_\infty
\bigl(
|p-q|+|q-r|+|r-p|
\bigr)
\\
&\le
2\|f\|_\infty
\bigl(
|p-r|+|q-r|
\bigr).
\end{aligned}
\]
Therefore,
\begin{equation}\label{eq:Ln_bound_corrected}
|\mathcal Q_n|
\le
2\|f\|_\infty
\left(
|\nabla_x\delta_{n-1}|
+
|\nabla_x\Phi_{n-1}|
\right),
\qquad n\ge1.
\end{equation}

\smallskip
\emph{Step 3: weighted energy estimate for $B_n$.}

We reverse time in~\eqref{eq:Bn_PDE_explicit}, so that the terminal
condition becomes a vanishing initial condition. Testing the resulting
equation with $B_n$ and integrating over $\mathbb R^d$, we obtain
\[
\begin{aligned}
\frac{\de}{\de t}\|B_n(t,\cdot)\|_2^2
+
\lambda\|\nabla_xB_n(t,\cdot)\|_2^2
\le{}&
2K_0
\int_{\mathbb R^d}
|\nabla_xB_n|\,|B_n|\,\de x
\\
&+
4\|f\|_\infty
\int_{\mathbb R^d}
|\nabla_x\Phi_{n-1}|\,|B_n|\,\de x
\\
&+
4\|f\|_\infty
\int_{\mathbb R^d}
|\nabla_x\delta_{n-1}|\,|B_n|\,\de x,
\end{aligned}
\]
where
\[
K_0
:=
\frac12
\|\operatorname{div}(\sigma\sigma^\top)\|_\infty
+
\|f\|_\infty.
\]
For any $\varepsilon>0$, Young's inequality gives
\[
2K_0ab
\le
\frac{\lambda}{4}a^2
+
\frac{4K_0^2}{\lambda}b^2
\]
and
\[
4\|f\|_\infty ab
\le
\varepsilon a^2
+
\frac{4\|f\|_\infty^2}{\varepsilon}b^2.
\]
Hence, for almost every $t\in[0,T]$,
\begin{equation}\label{eq:Bn_diff_ineq}
\begin{split}
\frac{\de}{\de t}\|B_n(t,\cdot)\|_2^2
+
\frac{3\lambda}{4}
\|\nabla_xB_n(t,\cdot)\|_2^2
\le{}&
K\|B_n(t,\cdot)\|_2^2
\\
&+
\varepsilon
\|\nabla_x\Phi_{n-1}(t,\cdot)\|_2^2
\\
&+
\varepsilon
\|\nabla_x\delta_{n-1}(t,\cdot)\|_2^2,
\end{split}
\end{equation}
where
\[
K
:=
\frac{4K_0^2}{\lambda}
+
\frac{8\|f\|_\infty^2}{\varepsilon}.
\]

For
\[
w\in
C([0,T];L^2(\mathbb R^d))
\cap
L^2(0,T;H^1(\mathbb R^d)),
\]
define
\[
\widetilde E(w)
:=
\sup_{t\in[0,T]}
e^{-Kt}\|w(t,\cdot)\|_2^2,
\qquad
\widetilde F(w)
:=
\int_0^T
e^{-Kt}\|\nabla_xw(t,\cdot)\|_2^2\,\de t.
\]
Multiplying~\eqref{eq:Bn_diff_ineq} by $e^{-Kt}$ and integrating in
time yields
\begin{equation}\label{eq:Bn_weighted}
\begin{split}
\widetilde E(B_n)
&\le
\varepsilon
\left(
\widetilde F(\Phi_{n-1})
+
\widetilde F(\delta_{n-1})
\right),
\\
\widetilde F(B_n)
&\le
\frac{4\varepsilon}{3\lambda}
\left(
\widetilde F(\Phi_{n-1})
+
\widetilde F(\delta_{n-1})
\right).
\end{split}
\end{equation}

The contraction estimate established in the proof of
Theorem~\ref{thm:PI_uniform} implies that there exists $C_\delta>0$
such that
\begin{equation}\label{eq:exact_increment_decay_practical}
\widetilde F(\delta_k)
\le
\int_0^T
\|\nabla_x\delta_k(t,\cdot)\|_2^2\,\de t
\le
C_\delta 2^{-k},
\qquad
k\ge0.
\end{equation}

\smallskip
\emph{Step 4: stability of the practical iterates.}

By~\eqref{eq:An_energy},
\[
\widetilde E(A_n)\le C_Ap_n^2,
\qquad
\widetilde F(A_n)\le\frac{C_A}{\lambda}p_n^2.
\]
Since $\Phi_n=A_n+B_n$, it follows from
\eqref{eq:Bn_weighted} that
\begin{equation}\label{eq:Phi_energy_contraction}
\begin{split}
\widetilde F(\Phi_n)
&\le
2\widetilde F(A_n)
+
2\widetilde F(B_n)
\\
&\le
\frac{2C_A}{\lambda}p_n^2
+
\frac{8\varepsilon}{3\lambda}
\left(
\widetilde F(\Phi_{n-1})
+
\widetilde F(\delta_{n-1})
\right).
\end{split}
\end{equation}
We now fix
\[
\varepsilon:=\frac{3\lambda}{32}.
\]
Then
\[
\widetilde F(\Phi_n)
\le
\frac{2C_A}{\lambda}p_n^2
+
\frac14\widetilde F(\Phi_{n-1})
+
\frac14\widetilde F(\delta_{n-1}).
\]
Because the two schemes have the same initial policy,
$\hat v_0=v_0$, and hence $\Phi_0=A_0$. Therefore,
\[
\widetilde F(\Phi_0)
\le
\frac{C_A}{\lambda}p_0^2.
\]
Using the monotonicity $p_k\le p_n$ for $k\le n$ and
\eqref{eq:exact_increment_decay_practical}, we may iterate the preceding
inequality to obtain
\[
\begin{aligned}
\widetilde F(\Phi_n)
&\le
C p_n^2
+
C_\delta
\sum_{k=1}^n
4^{-(n-k+1)}2^{-(k-1)}
\\
&\le
C\left(p_n^2+2^{-n}\right).
\end{aligned}
\]
Returning to the first inequality in~\eqref{eq:Bn_weighted}, we similarly
obtain
\[
\begin{aligned}
\widetilde E(\Phi_n)
&\le
2\widetilde E(A_n)
+
2\widetilde E(B_n)
\\
&\le
2C_Ap_n^2
+
2\varepsilon
\left(
\widetilde F(\Phi_{n-1})
+
\widetilde F(\delta_{n-1})
\right)
\\
&\le
C\left(p_n^2+2^{-n}\right).
\end{aligned}
\]
Removing the exponential weight therefore gives
\begin{equation}\label{eq:Phi_global}
\sup_{t\in[0,T]}
\|\Phi_n(t,\cdot)\|_2
\le
C\left(p_n+2^{-n/2}\right).
\end{equation}

\smallskip
\emph{Step 5: conclusion.}

By Theorem~\ref{thm:PI_uniform}, the exact policy-iteration sequence
satisfies
\[
\sup_{t\in[0,T]}
\|v_n(t,\cdot)-v(t,\cdot)\|_2
\le
C\rho^n,
\qquad
\rho=\frac1{\sqrt2}.
\]
Since
\[
\tilde v_n-v
=
\Phi_n+(v_n-v)
\]
and $2^{-n/2}=\rho^n$, estimate~\eqref{eq:Phi_global} yields
\[
\sup_{t\in[0,T]}
\|\tilde v_n(t,\cdot)-v(t,\cdot)\|_2
\le
C\bigl(p_n+\rho^n\bigr).
\]
This completes the proof.
\end{proof}

\begin{rem}[Behavior of the constants as $\lambda\to0$]\label{rem:small_lambda}
The guarantees of Theorems~\ref{thm:PI_uniform} and~\ref{thm:global_error} hold for every fixed $\lambda>0$, and no positive lower bound on $\lambda$ is required for their validity. The constants, however, are not uniform in $\lambda$. With the Young parameters chosen proportionally to $\lambda$ (e.g., $\varepsilon=3\lambda/8$ in the proof of Theorem~\ref{thm:PI_uniform}), the contraction rate $\rho=1/\sqrt2$ is independent of $\lambda$, but the Gronwall exponents scale as $K=O(1/\lambda)$, so the prefactor $C$ grows like $e^{cT/\lambda}$ as $\lambda\to0$. The estimates therefore degenerate in the vanishing-viscosity limit and do not extend to degenerate or first-order HJI equations; this limitation motivates the small-diffusion experiments in Section~\ref{sec:additional_experiments} and the discussion in Section~\ref{sec:conclusion}.
\end{rem}

This result highlights a key advantage of the proposed framework: unlike black-box direct PINN solvers, the iterative structure permits explicit $L^2$-error estimates with provable rates. This is particularly important in nonconvex settings, where bounding the solution error is otherwise analytically intractable.

Furthermore, each iteration returns a practical, near-optimal feedback policy via simple pointwise minimax updates, eliminating the need for any additional optimization.

\section{Experimental results}\label{sec:exp}

To demonstrate the effectiveness of our policy iteration scheme for solving nonconvex viscous HJI equations, we consider a two-dimensional optimal path planning problem in the presence of a moving obstacle. The setting is cast as a two-player zero-sum differential game, where the robot aims to reach a target while minimizing cost, and the environment acts as an adversarial player introducing worst-case disturbances. To assess scalability, we also apply our method to a high-dimensional publisher-subscriber game, where multiple agents interact under stochastic dynamics with anisotropic noise. Although no ground truth is available in this setting, the learned value functions exhibit symmetric structures that align with the problem's design.

\subsection{Implementation setup}\label{setup}

This section details the implementation aspects of our work, including the neural network architecture, training configuration, policy iteration setup, and construction of reference solutions.

\paragraph{Neural network architecture.}
As demonstrated in~Section~\ref{subsec:pi}, we parameterize the value function $v_n(t, x; \theta_n) = g(x) + (T - t) \mathcal{N}_n(t, x; \theta_n)$, thereby the terminal condition is satisfied automatically. We model $\mathcal{N}_n$ using a fully connected feedforward neural network (see Appendix~\ref{sec:impl_details} for details). In the neural network, we employ sinusoidal activation functions:
\[
\phi_i(v) = \sin(W_i v + b_i),
\]
as they are known to capture high-frequency structure and gradients more effectively than traditional activations~\cite{sitzmann2020implicit}. In high-dimensional reachability tasks, sine-based networks have also been shown to reduce mean squared errors by an order of magnitude compared to ReLU or tanh activations~\cite{bansal2021deepreach}.

\paragraph{Training configuration.}

We use the Adam optimizer~\cite{kingma2014adam} with a fixed learning rate of $0.001$. All implementations are carried out using the JAX framework, which enables efficient automatic differentiation and GPU acceleration. All experiments were conducted on a workstation equipped with an AMD Ryzen Threadripper PRO 9965WX CPU and an NVIDIA RTX PRO 6000 Blackwell Max-Q GPU.

\paragraph{Policy iteration setup.}
Our policy iteration scheme follows Algorithm~\ref{algo:main}. Problem-specific settings for the number of epochs $E$, policy updates $M$, and collocation points are detailed in Appendix~\ref{sec:impl_details}.

\paragraph{Reference solution.}
To quantitatively evaluate the accuracy of the learned value functions, we compute reference solutions to the HJI equation using an explicit finite difference method (FDM) on a uniform grid. We employ backward time integration with central differences in space, and discretize the diffusion term using a second-order scheme, with homogeneous Neumann boundary conditions imposed on the spatial domain. To maintain stability in the presence of diffusion, the time and space steps are chosen such that $\Delta t = (\Delta x)^2$. The resulting numerical solutions approximate the viscosity solution and are used to compute the relative $L^2$-errors.

\paragraph{Evaluation notations.}
To facilitate the comparative analysis using the aforementioned error metrics, we define the following notations for the different value function representations:
\begin{enumerate}[label=(\roman*), leftmargin=*]
    \item $v^{PI}_M(t, x; \theta_M)$: the value function obtained from our policy iteration-based PINN after $M$ policy updates.
    \item $v^D(t, x; \theta)$: the value function obtained from the direct PINN approach.
    \item $v^R(t, x)$: the reference solution computed using the FDM.
\end{enumerate}

\subsection{Two-dimensional optimal path planning with a moving obstacle}\label{sec:pathplanning}
Let \( X(s) \in \mathbb{R}^2 \) denote the position of the robot at time \( s \in [t, T] \). The system dynamics follow a controlled stochastic differential equation:
\[
\begin{cases}
    \begin{aligned}
    	\mathrm{d}X(s) &= (a(s)+b(s))\mathrm{d}s+ \sigma \mathrm{d} W_s \quad\text{for}\quad s \in (t,T],\\
    	X(t)&=x \in \mathbb{R}^2,
    \end{aligned}
\end{cases}
\]
where $a(s)$ satisfying $|a(s)| \leq 1$ denotes the control input of the robot (Player I), $b(s)$ satisfying $|b(s)|\leq \delta$ for some $\delta>0$ represents the adversarial disturbance (Player II), and \( \sigma \in \mathbb{R}^{2 \times 2} \) is the noise matrix. Given $x_{\text{goal}}\in \mathbb{R}^2$ and weights $\lambda_i$ for $i=1,2,3$, the cost functional is defined as:
\[
J(t,x;a, b) = \mathbb{E} \left [ \int_t^T  ( \lambda_1 |a(s)|^2 + \lambda_2 \phi(s, X(s))  ) \mathrm{d}s + \lambda_3 |X(T) - x_{\mathrm{goal}}|^2 \right ],
\]
where the obstacle penalty function \( \phi(t, x) \) is given by:
\[
\phi(s, x) = \exp \left( -\frac{|x - x_{\mathrm{obs}}(s)|^2}{2\varepsilon^2}  \right),
\quad x_{\mathrm{obs}}(s) = \begin{bmatrix} 0.5 \cos(\pi s) \\ 0.5 \sin(\pi s) \end{bmatrix}.
\]

We define the value function
\[
v(t,x) = \sup_{\beta \in \Gamma_t} \inf_{a \in \mathcal{A}_t} \mathbb{E} \left [ \int_t^T c(s, x(s), a(s), \beta[a](s)) \mathrm{d}s + g(x(T)) \right ],
\]
which is known to satisfy the viscous HJI equation
\[
\begin{cases}
	\partial_t v + H(t,x,\nabla_x v) = - \tfrac{1}{2} \mathrm{Tr}(\sigma \sigma^\top D_{xx}^2 v) \quad &\text{in}\quad [0,T)\times\mathbb{R}^d \\
	v(T,x)=g(x) \quad&\text{on}\quad \mathbb{R}^{d},
\end{cases}
\]
with
\[
H(t,x,p) = \sup_{\textbf{b} \in B} \inf_{\textbf{a} \in A}  [ \lambda_1 |\textbf{a}|^2 + \lambda_2 \phi(t,x) + p\cdot (\textbf{a} + \textbf{b})].
\]
for some $\lambda_1,  \lambda_2>0$. Here, the Hamiltonian takes the closed form given by
\[
H(t,x,p) =
\begin{cases}
	- \dfrac{1}{4\lambda_1} |p|^2 + \lambda_2 \phi(t,x) + \delta |p|, \quad&\text{if}\quad |p| \leq 2\lambda_1, \\
	- |p| + \lambda_1 + \lambda_2 \phi(t,x) + \delta |p|, \quad&\text{otherwise},
\end{cases}
\]
which is nonconvex as illustrated in Appendix~\ref{app:nonconvex_example}.

We solve the HJI equation using the PINN trained to minimize the residual of the PDE. Once the value function $ v(t, x; \theta) $ is learned, the optimal control policy is recovered via:
\[
a^*(t,x) =
\begin{cases}
	- \dfrac{1}{2\lambda_1} \nabla_x v(t, x; \theta), \quad& \text{if}\quad \left| \dfrac{1}{2\lambda_1} \nabla_x v(t, x; \theta) \right| \leq 1, \\
	- \dfrac{\nabla_x v(t, x; \theta)}{|\nabla_x v(t, x; \theta)|},\quad & \text{otherwise}.
\end{cases}
\]

We set the simulation domain to \( x \in [-1,1]^2 \) and the terminal time to \( T = 1.0 \). The optimization is configured with \( \delta = 0.1 \), \( \lambda_1 = 0.1 \), \( \lambda_2 = 100 \), \( \lambda_3 = 10 \), \( \varepsilon = 0.3 \), and the diffusion matrix is given by \( \sigma = 0.1 I_2 \). The target position is fixed at \( x_{\mathrm{goal}} = (0.9, 0.9) \).

\paragraph{Learning result.}
In the 2D moving obstacle example, the finite difference solution is computed on an extended spatial domain to reduce boundary artifacts. The restricted solution over the target region is then used as a quantitative baseline for evaluating the learned value function.

\begin{figure}[ht!]
	\centering
	\includegraphics[width=\textwidth]{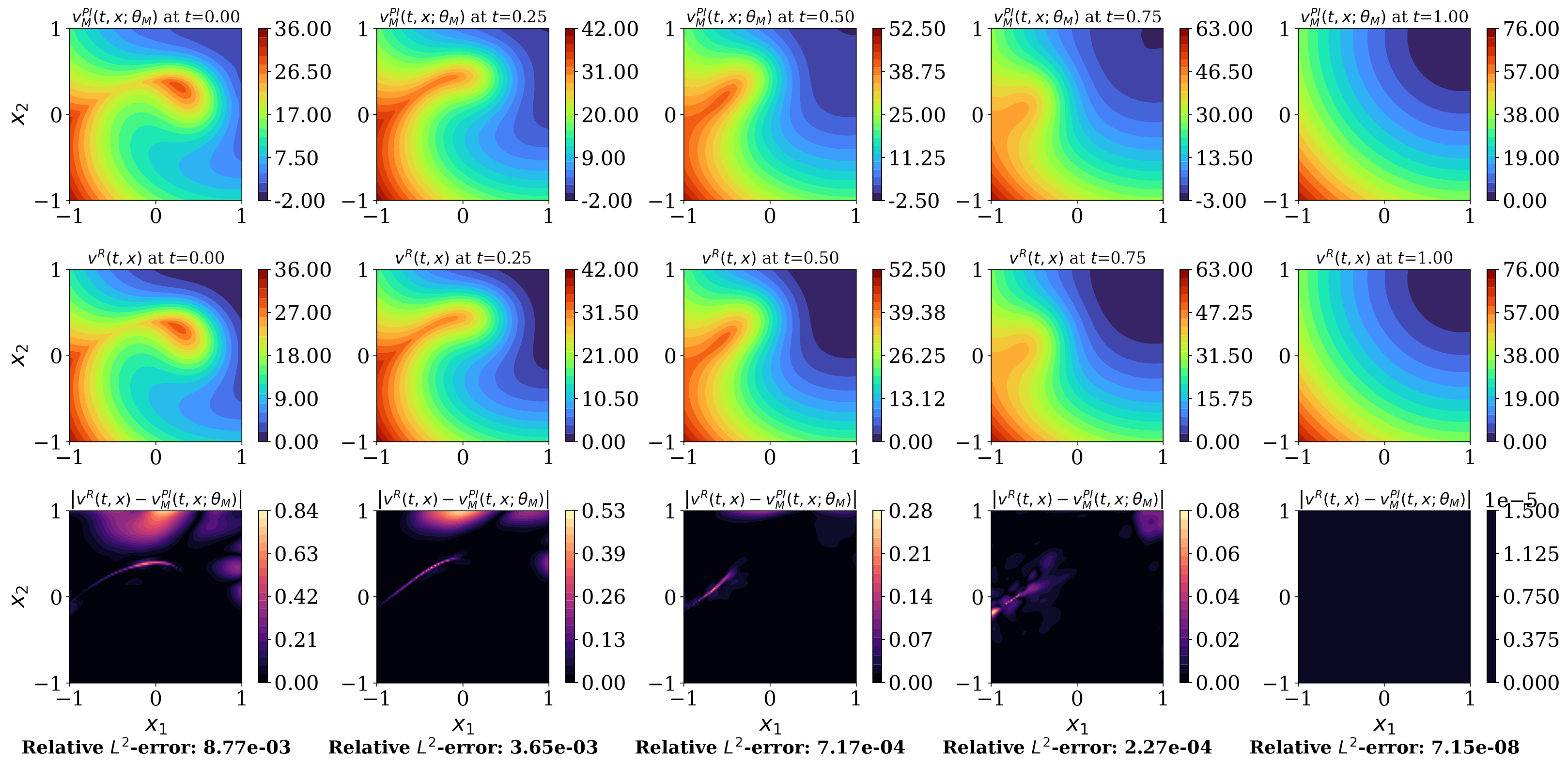}
	\caption{ 
		Comparison of the policy-iterative PINN solution and the reference FDM solution for the two-dimensional optimal path planning problem with a moving obstacle. Shown are the predicted value functions at selected times $t = 0.00$, $0.25$, $0.50$, $0.75$, and $1.00$, along with absolute error plots and corresponding relative $L^2$-error metrics.
	}
	\label{fig:opp:pinn_fdm_comparison}
\end{figure}

Figure~\ref{fig:opp:pinn_fdm_comparison} compares the learned value function from the policy-iterative PINN with the reference solution obtained via finite differences, at five representative time points. The first two rows show the predicted and reference value functions, respectively, while the bottom row displays the pointwise absolute error, along with the corresponding relative $L^2$-error. At all evaluated time instances, the policy-iterative PINN achieves consistently low errors, with relative $L^2$-errors on the order of $10^{-3}$ or lower, demonstrating strong agreement with the reference solution.

\begin{figure}[ht!]
	\centering
	\includegraphics[width=\textwidth]{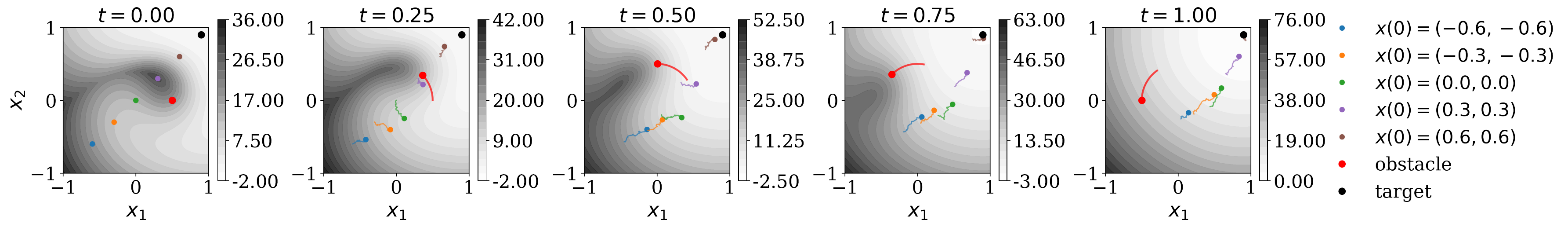}
	\caption{
		Time evolution of optimal trajectories derived from the policy-iterative PINN solution. Robots initialized at various positions navigate toward the target while avoiding a moving obstacle.
	}
	\label{fig:opp:traj_snapshots}
\end{figure}

Figure~\ref{fig:opp:traj_snapshots} illustrates the time evolution of optimal trajectories computed from the learned value function obtained via the policy-iterative PINN. At each time step, control inputs are derived from the gradient of the learned value function, and the system dynamics are integrated using the Euler--Maruyama method~\cite{higham2001algorithmic}. The resulting trajectories show how the agents, starting from different initial positions, successfully avoid the moving obstacle and reach the target by the cost landscape encoded in the value function.

\subsection{Publisher-subscriber differential game}\label{section:PS}
Let \(X(s)\in\mathbb{R}^{2N+1}\) denote the system state at time \(s\in[t,T]\). We write
\[
X(s)=\left(X_0(s),X_1(s),\ldots,X_{2N}(s)\right),
\]
where $X_0(s)$ is the publisher state and the subscribers are grouped into $N$ two-dimensional pairs: the $i$th pair is $(X_{2i-1}(s),X_{2i}(s))$ for $i=1,\ldots,N$.
Thus the total dimension is $2N+1$. A central agent (publisher) influences many followers (subscribers), each aiming to stay close to the leader while being disturbed. The dynamics are given by the controlled SDE
\[
\begin{cases}
\mathrm{d}X(s)=f\left(X(s),u(s),d(s)\right)\mathrm{d}s+\sigma\mathrm{d}W_s,
\quad \text{for}\quad s\in(t,T],\\
X(t)=x\in\mathbb{R}^{2N+1},
\end{cases}
\]
where \(W_s\) is a standard \((2N+1)\)-dimensional Brownian motion and \(\sigma\in\mathbb{R}^{(2N+1)\times(2N+1)}\) is the diffusion matrix.

\paragraph{Control and disturbance.}
The control input $u(t)\in\mathcal U$ is chosen by Player~I and acts only on the subscriber coordinates, while the disturbance $d(t)\in\mathcal D$ is chosen by Player~II:
\[
\mathcal U:=\{\mathbf u\in\mathbb{R}^{2N}:\|\mathbf u\|_\infty\le 1\},\qquad
\mathcal D:=\{\mathbf d\in\mathbb{R}^{2N}:\|\mathbf d\|_\infty\le 1\}.
\]

\paragraph{Drift structure.}
The drift is written in the compact form
\[
f(x,u,d)=Ax+Bu+Cd+\psi(x),
\]
where
\[
A = e_1 e_1^\top - \mathbf{1}_{2N+1} e_1^\top + a I_{2N+1},\qquad
B = \begin{bmatrix} \mathbf 0_{1\times 2N} \\ b I_{2N} \end{bmatrix}.
\]
Here $e_1\in\mathbb{R}^{2N+1}$ denotes the standard vector, i.e., $e_1=[1,0,\ldots,0]^\top$, $\mathbf{1}_{2N+1}\in\mathbb{R}^{2N+1}$ denotes the all-ones vector, i.e., $\mathbf{1}_{2N+1}=[1,\ldots,1]^\top$ and $\mathbf 0_{1\times 2N}$ denotes the zero row vector indicating that the control input does not act on the publisher coordinate $x_0$ and affects only the subscriber coordinates. 

\paragraph{Nonlinear interaction.}
To introduce genuine state-dependent nonlinearity while preserving the block structure, we include the interaction term
\[
\psi(x)=
\begin{bmatrix}
\alpha\sin(x_0)\\
-\beta x_0 \mathbf 1_{2N}
\end{bmatrix}\circ (x\circ x).
\]
This term induces nonlinear coupling from the publisher to each subscriber pair,
while avoiding cross-coupling between different subscriber blocks.

\paragraph{Rotated disturbance and nonconvex Isaacs structure.}
The disturbance enters through a rotated block structure,
\[
C=\begin{bmatrix}
\mathbf 0_{1\times 2N} \\ c \mathrm{diag}(R(\phi_1),\ldots,R(\phi_N))
\end{bmatrix} \in \R^{(2N+1)\times 2N},
\qquad
R(\phi)=
\begin{bmatrix}
\cos\phi & -\sin\phi\\
\sin\phi & \cos\phi
\end{bmatrix}.
\]
This rotation breaks axis alignment between control and disturbance channels.
As a consequence, the Isaacs term in the Hamiltonian becomes generally neither convex nor concave
in the gradient variable $p$; see Appendix~\ref{app:nonconvex_example}.

\paragraph{Terminal cost and value function decomposition.}
The terminal cost is separable across publisher-subscriber pairs:
\[
g(x)=\sum_{i=1}^{N} g_i(P_i x),
\qquad
g_i(z):=\tfrac12\left(z_1^2+z_2^2+z_3^2-r^2\right),\quad z\in\mathbb{R}^3.
\]
Here \(P_i:\mathbb{R}^{2N+1}\to\mathbb{R}^3\) is the projection onto the \((x_0,x_{2i-1},x_{2i})\) coordinates,
\[
P_i x := \left[x_0,  x_{2i-1},  x_{2i}\right]^\top,
\]
so that \(z=P_i x\) and \(g_i(P_i x)=\tfrac12\left(x_0^2+x_{2i-1}^2+x_{2i}^2-r^2\right)\). Together with the unidirectional coupling and the block-diagonal diffusion
$$
\sigma=\mathrm{diag}\left(\sigma_0,\Sigma_1,\ldots,\Sigma_N\right)\in \mathbb{R}^{(2N+1)\times (2N+1)},
$$
where $\sigma_0\in\mathbb{R}$ and each $\Sigma_i\in\mathbb{R}^{2\times2}$ (e.g., $\sigma_0=0.1$ and $\Sigma_i=0.1 I_2$). This structure yields an exact decomposition of the value function,
\[
v(t,x)=\sum_{i=1}^{N} v_i(t,x_0,x_{2i-1},x_{2i}),
\]
where each $v_i$ solves a three-dimensional HJI equation over the $(x_0,x_{2i-1},x_{2i})$ subspace; see Appendix~\ref{proof:value_decomp} for the proof.

\emph{Hamiltonian structure.} 
The corresponding Hamiltonian takes the form of
\[
H(x, p) = p^\top (A x + \psi(x)) - \|B^\top p\|_1 + \|C^\top p\|_1,
\]
where $\|z\|_1:=\sum_k|z_k|$. Note that, due to the rotated disturbance map $C$, the Isaacs term
$-\|B^\top p\|_1+\|C^\top p\|_1$ is generally neither convex nor concave in $p$;
see Figure~\ref{fig:hamiltonian_nonconvex_app}(b) in  Appendix~\ref{app:nonconvex_example} for the case $N=1$ . 

The HJI equation is solved using a PINN, and the optimal control and disturbance policies are recovered from the gradient $p=\nabla_x v(t,x)$ of the learned value function as
$$
u^*(x)=-\mathrm{sign}(B^\top p),\qquad
d^*(x)=\mathrm{sign}(C^\top p),
$$
where the sign function is applied componentwise.

To solve the HJI equation, we consider both the proposed policy-iterative PINN and a direct PINN baseline. Both models share the same architecture, but differ in loss formulation and training schedule (see Appendix~\ref{sec:impl_details}). The policy-iterative PINN minimizes the residual under fixed policies (see Section~\ref{setup}), while the direct PINN substitutes the closed-form Hamiltonian into the objective:
\begin{align*}
	&\mathcal{L}_{\text{Direct}}(\theta) =
	\\&\quad\frac{1}{N_{\mathrm{int}}} \sum_{j=1}^{N_{\mathrm{int}}}
	[
	\partial_t v(t_j, x_j;\theta)
	+ H(x_j, \nabla_x v(t_j, x_j; \theta))
	+ \frac{1}{2} \operatorname{Tr} ( \sigma \sigma^\top D^2_{xx} v(t_j, x_j; \theta)  )
	]^2.
\end{align*}
Although the direct approach avoids explicit policy updates, the resulting loss involves non-smooth and nonlinear terms such as $\ell_1$ norms of gradients, which may result in instability during training. The policy-iterative formulation mitigates this by decoupling the optimization into simpler fixed-policy subproblems. 

In all problem settings described below, we assess accuracy by comparing to a reference solution defined over an extended domain and restricted to the target region. Training is performed on the same domain to improve boundary behavior. Additional results on error convergence and training time are provided in Section~\ref{sec:additional_results}.

We set the simulation domain to $x\in[-0.5,0.5]^{2N+1}$ with terminal time $T=0.5$ and parameters
$a=1$, $b=1$, $c=\tfrac12$, $\alpha=-2$, $\beta=2$, and $\phi_i=\pi/4$. The diffusion matrix is defined as \( \sigma = 0.1  I + P_\epsilon \), where \(P_\epsilon\) is a symmetric matrix with zero diagonal and off-diagonal entries sampled from \( \mathcal{U}(0, \epsilon) \). This formulation yields an isotropic setting when $\epsilon = 0$ and becomes anisotropic otherwise.

\subsubsection{3D example - isotropic noise}

We first consider the isotropic three-dimensional case ($\epsilon=0$), corresponding to the simplest ($N=1$) instance of our model, and use it as an initial testbed before moving to higher-dimensional settings. For visualization, we report 2D contour plots on the diagonal slice $(x_0, s)$ with $x_1=x_2=s$.

\paragraph{Learning result}

\begin{figure}[ht!]
	\centering
	\includegraphics[width=\textwidth]{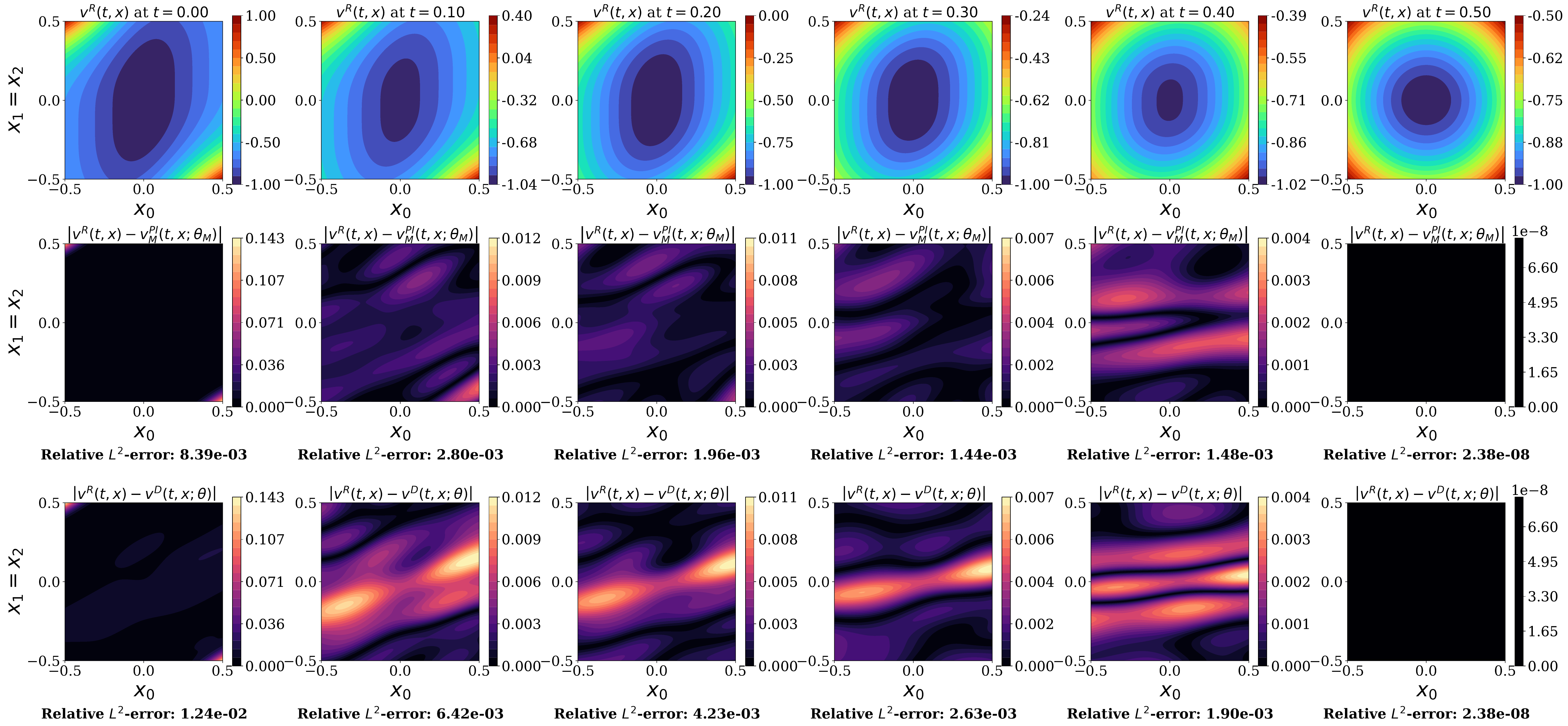}
	\caption{
    Comparison of policy-iterative and direct PINN approaches on the 3D isotropic benchmark problem.
    Columns show 2D slices at $t=0.0, 0.1, \ldots, 0.5$ (with $x_1=x_2$); rows (top to bottom) show the reference value function and the absolute error maps for the policy-iterative and direct PINN methods.
    Relative $L^2$-errors are reported below each panel; all panels correspond to the best-performing seed.
    }
	\label{fig:error:3d:isotropic}
\end{figure}

\begin{table}[ht!]
	\caption{Relative $L^2$-errors over time in the 3D isotropic setting for policy-iterative and direct PINN methods. Results are summarized by the mean, standard deviation, and best-seed (minimum) error over five random seeds.
		These results correspond to the time slices shown in Figure~\ref{fig:error:3d:isotropic}.}
	\label{tab:error:3d:isotropic}
	\begin{center}
		\vskip-.3truecm
        \resizebox{\textwidth}{!}{%
		\begin{tabular}{cccccccc}
			\hline\hline
			Method & Metric & $t=0.0$ & $t=0.1$ & $t=0.2$ & $t=0.3$ & $t=0.4$ & $t=0.5$ \\
			\noalign{\hrule height 1.5pt}
			\multirow{3}{*}{PI} 
			& Mean        & 1.074e-02 & 5.209e-03 & 3.837e-03 & 2.532e-03 & 1.599e-03 & 2.337e-08 \\
			& Std         & (2.45e-03) & (2.48e-03) & (1.72e-03) & (9.46e-04) & (2.86e-04) & (0.00e-00) \\
%			\cline{2-8}
			& \bf{Best}        & \bf{8.393e-03} & \bf{2.800e-03} & \bf{1.963e-03} & \bf{1.436e-03} & \bf{1.483e-03} &\bf{2.337e-08} \\
			\noalign{\hrule height 1pt}
			\multirow{3}{*}{Direct} 
			& Mean        & 1.431e-02 & 7.869e-03 & 5.392e-03 & 3.261e-03 & 1.979e-03 & 2.337e-08 \\
			& Std         & (1.77e-03) & (1.66e-03) & (9.36e-04) & (4.87e-04) & (1.25e-04) & (0.00e-00) \\
%			\cline{2-8}
			& \bf{Best}        & \bf{1.244e-02} & \bf{6.418e-03} & \bf{4.233e-03} & \bf{2.626e-03} & \bf{1.900e-03} & \bf{2.337e-08} \\
			\noalign{\hrule height 1.5pt}
		\end{tabular}
        }
	\end{center}
	\vskip-.2truecm
\end{table}

Figure~\ref{fig:error:3d:isotropic} shows the reference solution alongside the point-wise absolute errors of the policy-iterative and direct PINN methods using the best-performing seed at selected times $t = 0.0$, $0.1$, $\dots$, $0.5$. To provide a more comprehensive evaluation, Table~\ref{tab:error:3d:isotropic} summarizes the relative $L^2$-errors for all five random seeds, reporting the mean, standard deviation, and best-seed (minimum) error. At $t = 0$, the policy-iterative PINN achieves a best-seed relative $L^2$-error on the order of $10^{-3}$, indicating close agreement with the reference solution. Across the considered time steps, the policy-iterative PINN yields slightly improved accuracy in the mean error.
Note that in Tables~\ref{tab:error:3d:isotropic}--\ref{tab:error:11d:isotropic} both methods attain the identical relative $L^2$-error of $2.337\times10^{-8}$ at the terminal time $t=0.5$: since the ansatz \eqref{eq:ansatz} enforces the terminal condition $v(T,\cdot)=g$ as a hard constraint for both networks, the approximation error at $t=T$ vanishes up to floating-point precision.

\subsubsection{5D example}

We evaluate the two PINN methods on a five-dimensional ($N=2$) variant of the differential game in the isotropic setting ($\epsilon = 0$) using error tables, and additionally examine the anisotropic case ($\epsilon > 0$) via policy-iterative PINN visualizations. All configurations follow the high-dimensional setup described in Appendix~\ref{sec:impl_details}. To construct a reference solution, we exploit the separable structure of the isotropic case by summing two three-dimensional FDM solutions. We visualize the five-dimensional value function using the same 2D diagonal slice representation, plotting it over $(x_0, s)$ with $x_1=x_2=x_3=x_4=s$, which preserves the paired-subscriber structure $(x_0,x_1,x_2)$ and $(x_0,x_3,x_4)$ while enabling 2D contour plots.

\paragraph{Learning result}

\begin{figure}[ht!]
	\centering
	\includegraphics[width=\textwidth]{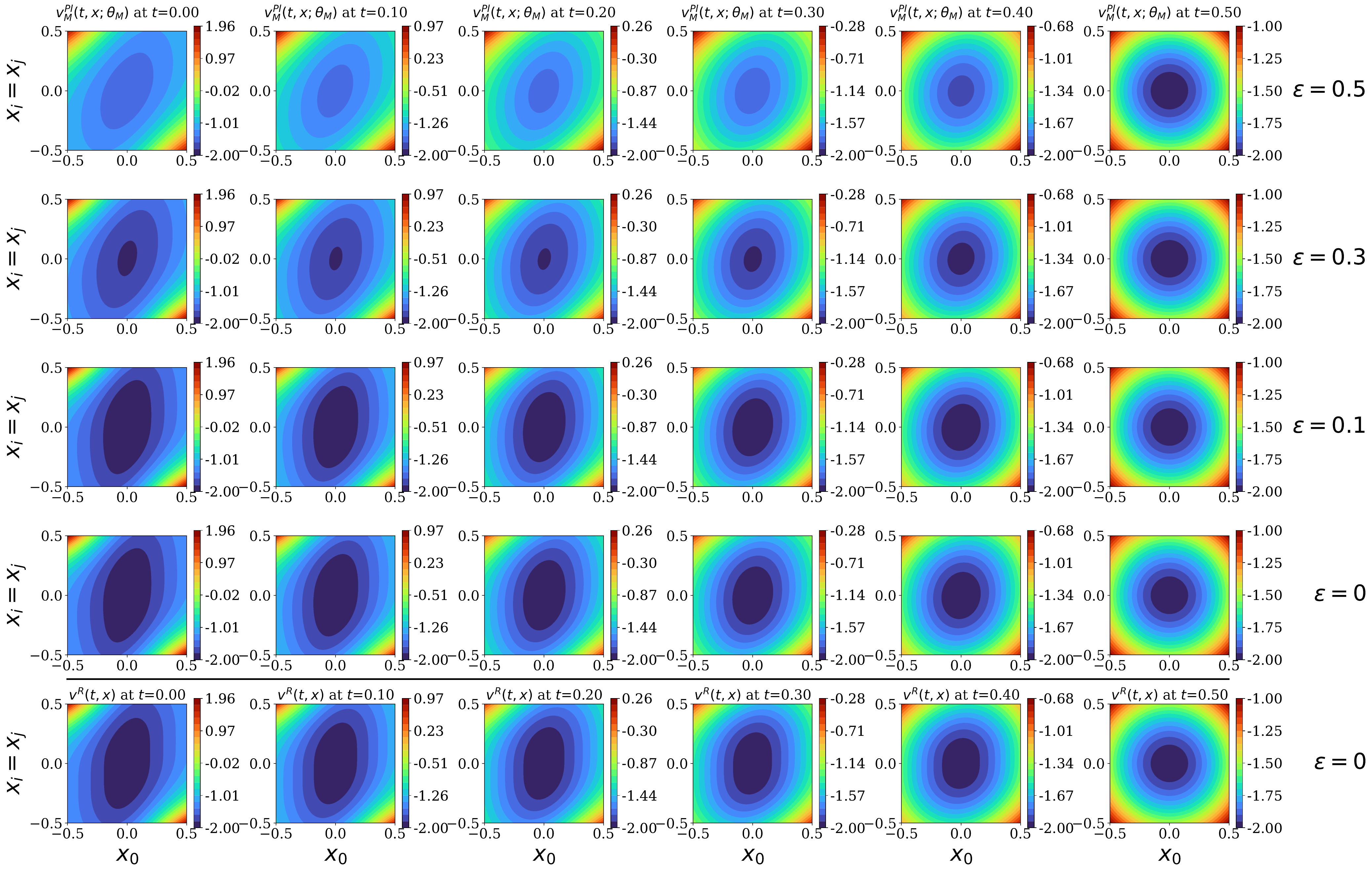}
	\caption{
		Policy-iterative PINN on the five-dimensional ($N=2$) anisotropic problem. Columns show the value functions at $t=0.0, 0.1,\ldots,0.5$. Rows (top to bottom) correspond to the learned $v^{PI}_M(t,x;\theta_M)$ for $\epsilon=0.5, 0.3, 0.1, 0.0$, while the last row shows the reference solution $v^R(t,x)$ for the isotropic case ($\epsilon=0.0$). As $\epsilon$ decreases, the learned contours progressively approach the isotropic reference.
	}
	\label{fig:error:5d:anisotropic}
\end{figure}

\begin{table}[ht!]
	\caption{Relative $L^2$-errors over time in the 5D isotropic setting for policy-iterative and direct PINN methods. Results are summarized by the mean, standard deviation, and best-seed (minimum) error over five random seeds.}
	\label{tab:error:5d:isotropic}
	\begin{center}
		\vskip-.3truecm
        \resizebox{\textwidth}{!}{%
		\begin{tabular}{cccccccc}
			\hline\hline
			Method & Metric & $t=0.0$ & $t=0.1$ & $t=0.2$ & $t=0.3$ & $t=0.4$ & $t=0.5$ \\
			\noalign{\hrule height 1.5pt}
			\multirow{3}{*}{PI} 
			& Mean        & 1.606e-02 & 8.295e-03 & 5.976e-03 & 5.254e-03 & 4.017e-03 & 2.337e-08 \\
			& Std         & (2.70e-03) & (2.29e-03) & (1.10e-03) & (4.45e-04) & (4.39e-04) & (0.00e-00) \\
%			\cline{2-8}
			& \bf{Best}        & \bf{1.274e-02} & \bf{5.388e-03} & \bf{5.006e-03} & \bf{5.553e-03} & \bf{4.111e-03} &\bf{2.337e-08} \\
			\noalign{\hrule height 1pt}
			\multirow{3}{*}{Direct} 
			& Mean        & 1.706e-02 & 7.683e-03 & 5.895e-03 & 5.608e-03 & 4.277e-03 & 2.337e-08 \\
			& Std         & (1.66e-03) & (1.52e-03) & (1.45e-03) & (1.14e-03) & (7.21e-04) & (0.00e-00) \\
%			\cline{2-8}
			& \bf{Best}        & \bf{1.542e-02} & \bf{5.865e-03} & \bf{3.995e-03} & \bf{3.717e-03} & \bf{2.928e-03} & \bf{2.337e-08} \\
			\noalign{\hrule height 1.5pt}
		\end{tabular}
        }
	\end{center}
	\vskip-.2truecm
\end{table}

Figure~\ref{fig:error:5d:anisotropic} reports the policy-iterative PINN value functions across time and varying levels of anisotropy on the 2D diagonal slice, plotting over $(x_0,s)$ with $x_1=x_2=x_3=x_4=s$, together with the summed isotropic reference on the same slice. In the isotropic case ($\epsilon=0$), where the summed reference is available via decomposition, we quantitatively compare the policy-iterative and direct PINN using the relative $L^2$-errors in Table~\ref{tab:error:5d:isotropic}, reported as the mean and standard deviation over five random seeds together with the best-seed (minimum) error across $t=0.0, 0.1,\ldots, 0.5$; overall, the policy-iterative PINN attains lower errors on average. For anisotropic settings ($\epsilon>0$), no ground-truth reference is available, so we assess performance qualitatively from the policy-iterative PINN contour evolution.

\subsubsection{11D example}

\begin{figure}[ht!]
	\centering
	\includegraphics[width=\textwidth]{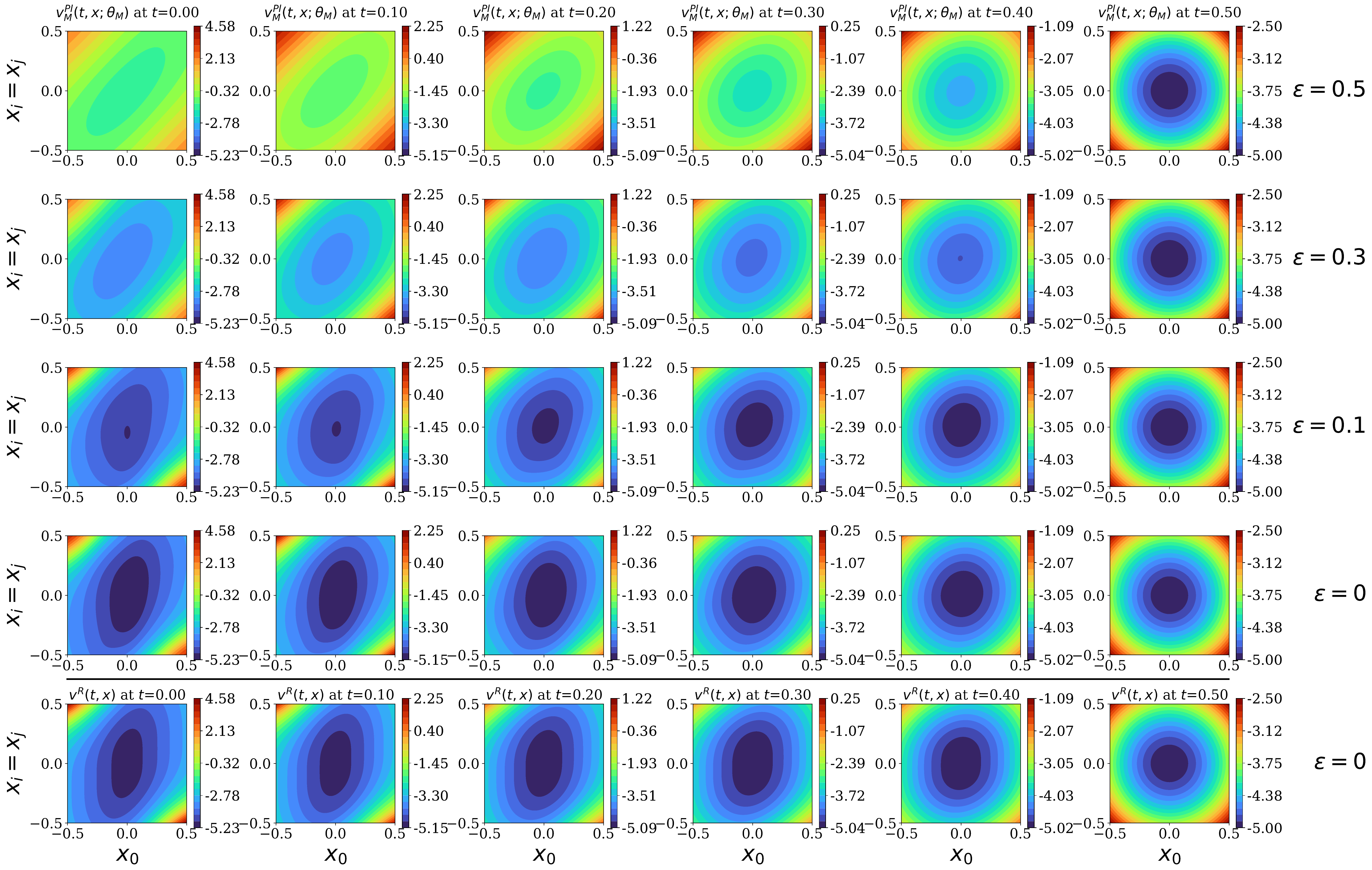}
\caption{Policy-iterative PINN on the eleven-dimensional ($N=5$) anisotropic problem. Columns show 2D diagonal-slice value functions at $t=0.0, 0.1, \ldots, 0.5$. Rows (top to bottom) correspond to the learned $v^{PI}_M(t,x;\theta_M)$ for $\epsilon = 0.5, 0.3, 0.1, 0.0$, while the last row shows the reference solution $v^R(t,x)$ for the isotropic case ($\epsilon = 0.0$), obtained via decomposition. All plots are shown on the slice $x_1 = x_2 = \cdots = x_{10} = s$ and visualized over $(x_0, s)$.
As $\epsilon$ decreases, the learned contours progressively approach the isotropic reference.
}

	\label{fig:error:11d:anisotropic}
\end{figure}

\begin{table}[ht!]
	\caption{Relative $L^2$-errors over time in the 11D isotropic setting for policy-iterative and direct PINN methods. Results are summarized by the mean, standard deviation, and best-seed (minimum) error over five random seeds.}
	\label{tab:error:11d:isotropic}
	\begin{center}
		\vskip-.3truecm
        \resizebox{\textwidth}{!}{%
		\begin{tabular}{cccccccc}
			\hline\hline
			Method & Metric & $t=0.0$ & $t=0.1$ & $t=0.2$ & $t=0.3$ & $t=0.4$ & $t=0.5$ \\
			\noalign{\hrule height 1.5pt}
			\multirow{3}{*}{PI} 
			& Mean        & 4.797e-02 & 3.062e-02 & 2.410e-02 & 1.902e-02 & 1.191e-02 & 2.337e-08 \\
			& Std         & (8.74e-03) & (4.25e-03) & (4.79e-03) & (4.91e-03) & (2.82e-03) & (0.00e-00) \\
%			\cline{2-8}
			& \bf{Best}        & \bf{3.556e-02} & \bf{2.525e-02} & \bf{2.028e-02} & \bf{1.700e-02} & \bf{1.309e-02} &\bf{2.337e-08} \\
			\noalign{\hrule height 1pt}
			\multirow{3}{*}{Direct} 
			& Mean        & 6.262e-02 & 3.527e-02 & 2.294e-02 & 1.705e-02 & 9.552e-03 & 2.337e-08 \\
			& Std         & (1.38e-02) & (2.35e-03) & (3.21e-03) & (2.12e-03) & (6.41e-04) & (0.00e-00) \\
%			\cline{2-8}
			& \bf{Best}        & \bf{4.152e-02} & \bf{3.103e-02} & \bf{2.622e-02} & \bf{1.789e-02} & \bf{8.581e-03} & \bf{2.337e-08} \\
			\noalign{\hrule height 1.5pt}
		\end{tabular}
        }
	\end{center}
	\vskip-.2truecm
\end{table}

We extend the evaluation to the eleven-dimensional ($N=5$) setting, under both isotropic ($\epsilon = 0$) and anisotropic ($\epsilon > 0$) noise, using the same experimental configuration and visualization strategy as in the five-dimensional case.

\paragraph{Learning result}

Figure~\ref{fig:error:11d:anisotropic} presents the policy-iterative PINN value functions in the 11D setting using the same visualization layout as in the 5D case, including the summed isotropic reference in the last row. In the isotropic case ($\epsilon=0$), we quantitatively compare the policy-iterative and direct PINN using the relative $L^2$-errors reported in Table~\ref{tab:error:11d:isotropic} (mean, standard deviation, and best-seed over five random seeds) across $t=0.0, 0.1,\ldots, 0.5$. As expected, the 11D problem exhibits slightly larger errors than the 5D case due to increased dimensionality. For anisotropic settings ($\epsilon>0$), no ground-truth reference is available, so we assess performance qualitatively from the policy-iterative PINN contour evolution.

\subsection{Further discussion on experimental results}\label{sec:additional_results}
While the preceding subsections report the final relative $L^2$-errors at selected time instances, they do not capture the optimization behavior of the two approaches during training. In this subsection, we therefore examine the training dynamics and computational efficiency of the proposed policy-iterative PINN, and compare them with those of the direct PINN baseline.

\paragraph{Optimization dynamics.}

\begin{figure}[t]
    \centering
    \begin{subfigure}[t]{0.32\textwidth}
        \centering
        \includegraphics[width=\linewidth]{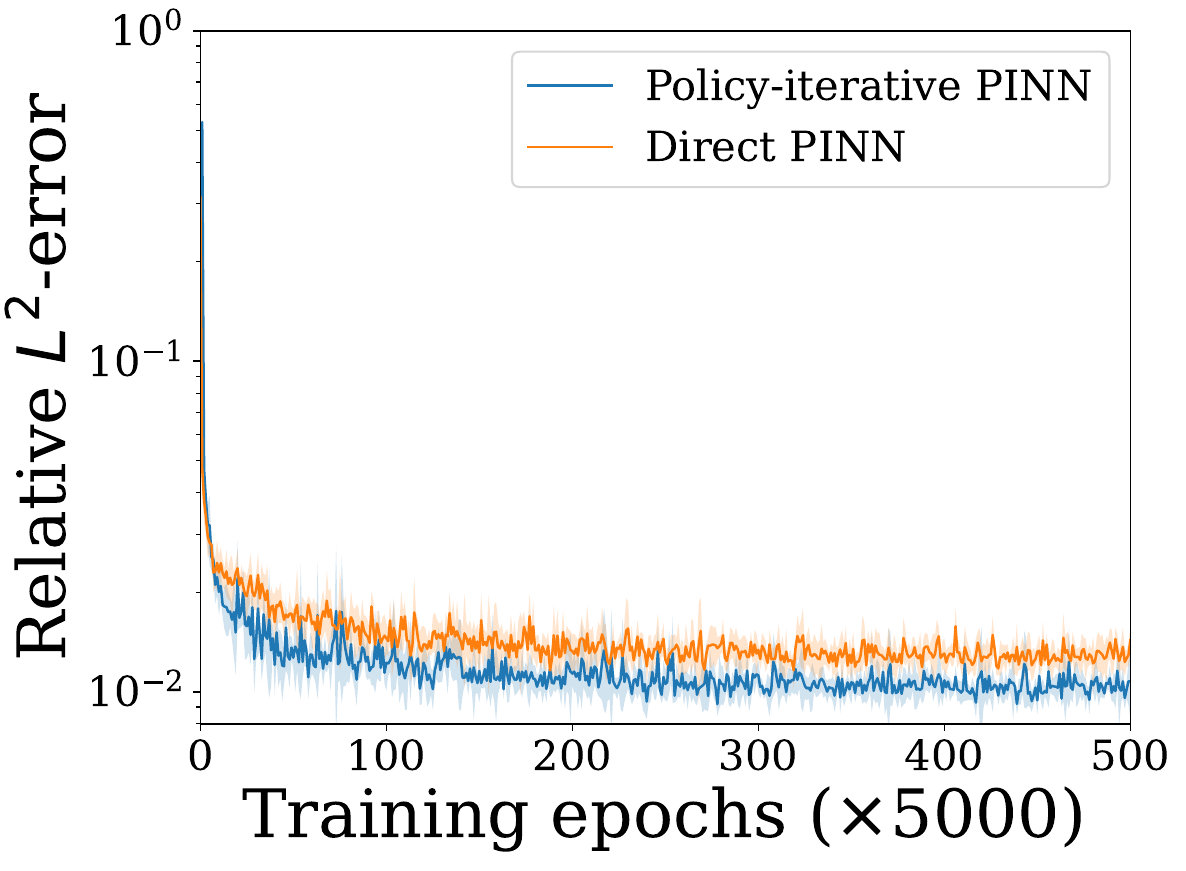}
        \subcaption{3D ($N=1$)}
        \label{fig:err_conv_3d}
    \end{subfigure}\hfill
    \begin{subfigure}[t]{0.32\textwidth}
        \centering
        \includegraphics[width=\linewidth]{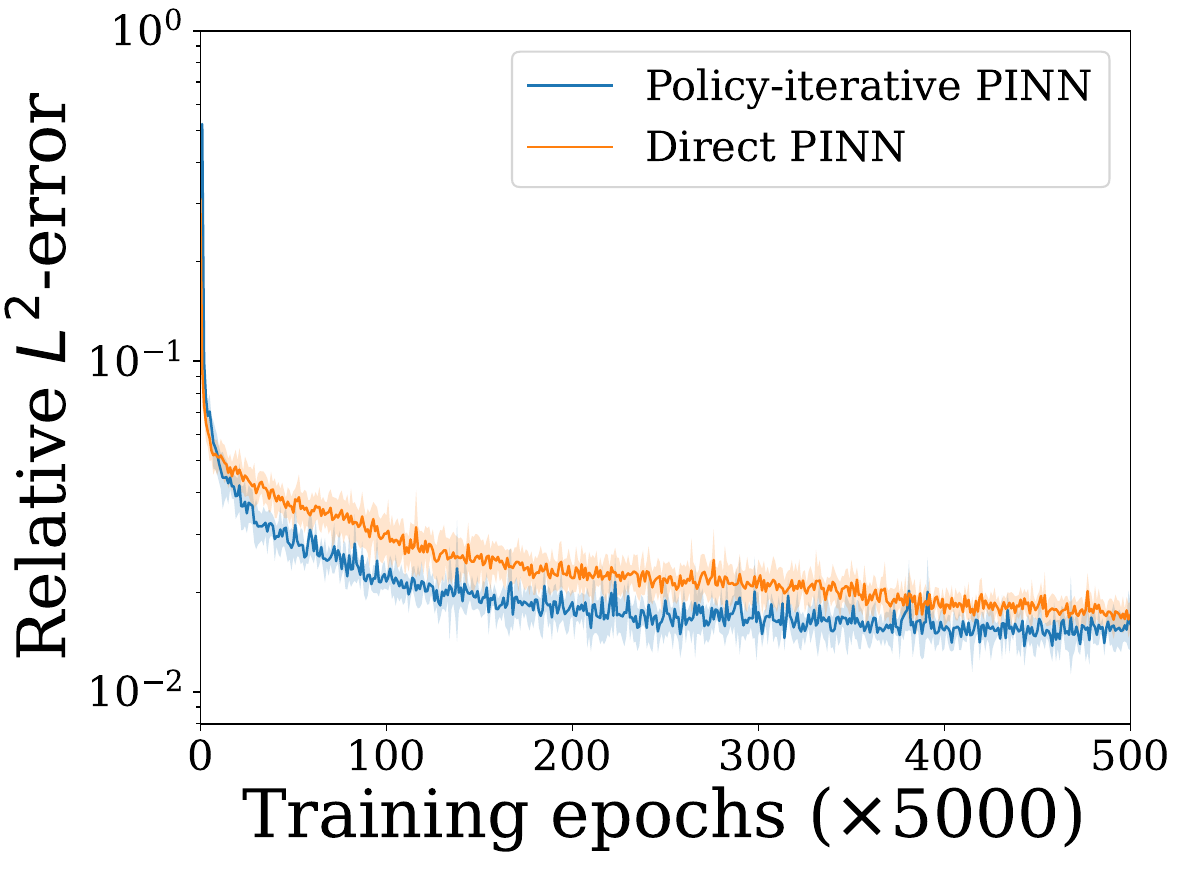}
        \subcaption{5D ($N=2$)}
        \label{fig:err_conv_5d}
    \end{subfigure}\hfill
    \begin{subfigure}[t]{0.32\textwidth}
        \centering
        \includegraphics[width=\linewidth]{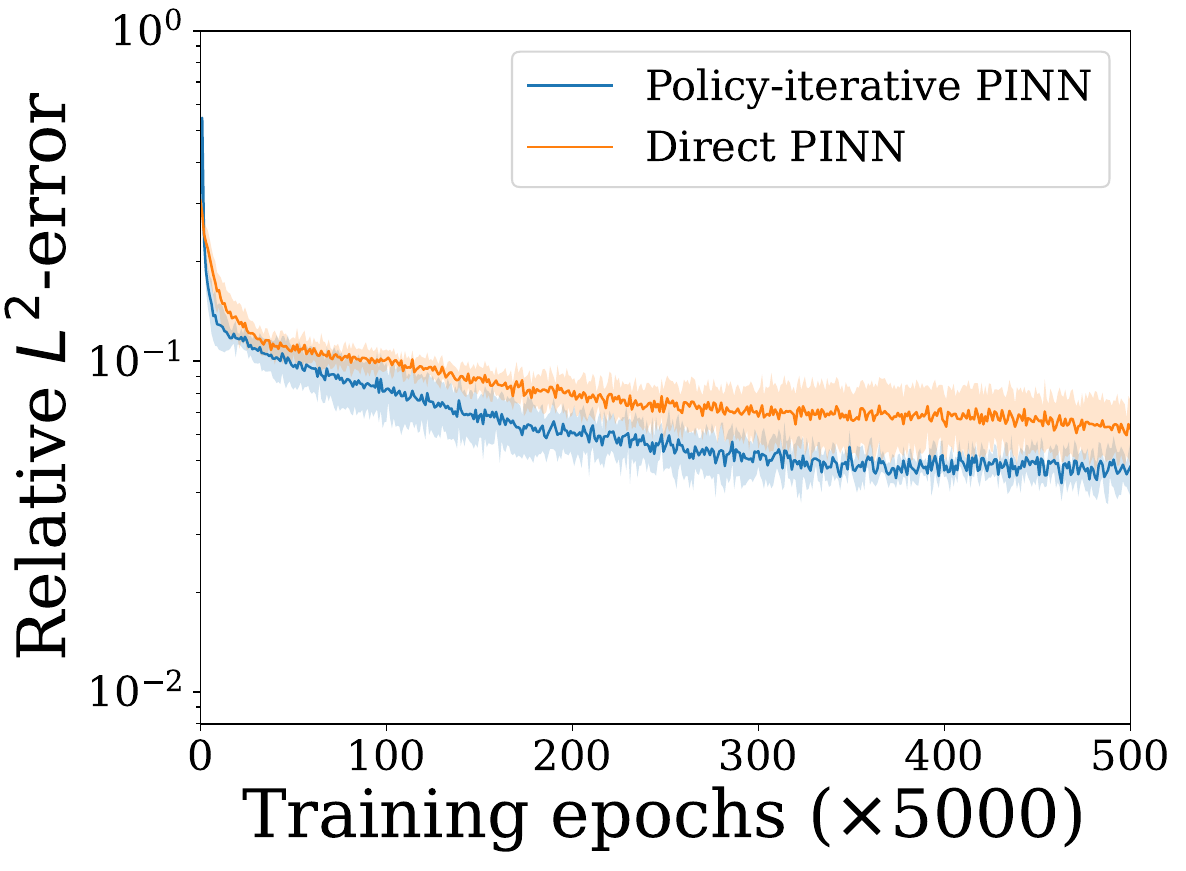}
        \subcaption{11D ($N=5$)}
        \label{fig:err_conv_11d}
    \end{subfigure}

    \caption{
    Relative $L^2$-error versus training progress for the policy-iterative PINN and the direct PINN in 3D ($N=1$), 5D ($N=2$), and 11D ($N=5$).
    Solid lines show the mean over five random seeds, and shaded regions indicate one standard deviation.
    }
    \label{fig:err_convergence}
\end{figure}

Figure~\ref{fig:err_convergence} shows the evolution of the relative $L^2$-error as a function of training progress for the 3D ($N=1$), 5D ($N=2$), and 11D ($N=5$) problems. For the policy-iterative PINN, the error is evaluated after each policy update, where the value network is trained for a fixed number of epochs with the policy frozen. For the direct PINN, the error is recorded at the same interval to ensure a fair comparison.

Across all tested dimensions, the policy-iterative PINN converges more rapidly and reaches a lower error floor than the direct PINN.
Moreover, the performance gap becomes more pronounced as the state dimension increases. In addition to lower mean errors, the policy-iterative PINN exhibits reduced variance across random seeds, indicating more stable optimization behavior.

\paragraph{Wall-clock efficiency.}

\begin{table}[ht!]
	\caption{Total training time(s) for policy-iterative and direct PINNs in the 3D, 5D, and 11D settings.}
	\label{tab:time:3d5d11d}
	\begin{center}
		\vskip-.3truecm
		\begin{tabular}{cccc}
			\hline\hline
			Method & 3D & 5D & 11D \\
			\noalign{\hrule height 1.5pt}
			  Policy-iterative PINN     & 2{,}819 & 9{,}221 & 83{,}761 \\
			\noalign{\hrule height 1pt}
			Direct PINN & 2{,}786 & 9{,}207 & 83{,}318 \\
			\noalign{\hrule height 1.5pt}
		\end{tabular}
	\end{center}
	\vskip-.2truecm
\end{table}

Importantly, this improvement in accuracy is not achieved at the expense of additional computational cost. Table~\ref{tab:time:3d5d11d} reports the total wall-clock training time for both methods in the 3D, 5D, and 11D settings. The training times are nearly identical across all dimensions,
demonstrating that the improved convergence of the policy-iterative PINN
is not due to increased compute, but rather to more favorable optimization conditioning.

The near-identical runtimes can be explained by examining where automatic differentiation is actually spent in the two methods. In the policy-iterative PINN, the policy-related computations consist of (i) the policy-improvement step, performed once per outer iteration after the value network has been trained under the current frozen policies, and (ii) the evaluation of the frozen feedback laws on refreshed collocation batches during inner training (see Section~\ref{subsec:pi}). Both computations require only the first-order state gradient $\nabla_x v$; they involve neither the second-order term $\operatorname{Tr}(\sigma\sigma^\top D^2_{xx}v)$ nor backpropagation of the PDE residual with respect to the network parameters. By contrast, the dominant cost in both methods is the residual-based training step itself, in which every optimizer step evaluates first- and second-order spatial derivatives of the network and backpropagates the residual through all trainable parameters. Since the two methods share the same architecture, number of optimizer steps, and collocation budget (Appendix~\ref{sec:impl_details}), the additional policy-related computations of the policy-iterative PINN are negligible relative to this dominant cost, and the wall-clock training times remain nearly identical.

\paragraph{Interpretation.}
This difference can be attributed to the structural distinction between the two formulations. The direct PINN enforces the Hamilton--Jacobi--Isaacs equation by embedding the nested minimax operator directly into the loss function. As a result, the optimization landscape is highly nonconvex and nonsmooth, notably due to the presence of $\ell_1$-type gradient terms and the implicit coupling between value approximation and policy optimization.
These effects become increasingly detrimental as the state dimension grows.

In contrast, the policy-iterative formulation decouples the Isaacs minimax structure into a sequence of fixed-policy subproblems.
Each policy evaluation step corresponds to solving a linear parabolic PDE with frozen controls, which yields a significantly smoother and better-conditioned loss surface. The subsequent policy improvement step updates the controls pointwise using gradients computed by automatic differentiation, without backpropagating through the minimax operator. This decoupling explains the smoother error decay and lower error floors
observed in Figure~\ref{fig:err_convergence}, and is consistent with the $L^2$-contractive behavior established in Section~\ref{sec:method}.

This favorable conditioning stems from a structural difference between the two loss functions. In our benchmark problems, the direct PINN does not solve an explicit inner policy-optimization problem during training, because the minimax Hamiltonian is available in closed form and is substituted directly into the residual. This closed-form Hamiltonian contains nonsmooth switching terms, such as absolute values of the value gradient, arising from the pointwise minimization and maximization; see, for instance, the piecewise Hamiltonian in Section~\ref{sec:pathplanning}. Every optimizer step of the direct PINN must differentiate the residual through these kinks. In the policy-iterative PINN, the control and disturbance policies are frozen during each value-update stage, so their contributions enter the residual linearly in the value gradient, and the loss minimized during inner training avoids the nonsmooth switching structure altogether. The advantage therefore lies in the structure of the optimization problem rather than in a reduction of computational work: the policy-iterative PINN achieves improved accuracy and stability under a computational budget comparable to that of the direct PINN, not by reducing wall-clock runtime.

\paragraph{Sensitivity to collocation density.}

\begin{figure}[h!]
    \centering
    %% 3D
    \begin{subfigure}[t]{0.49\textwidth}
        \centering
        \includegraphics[width=\textwidth]{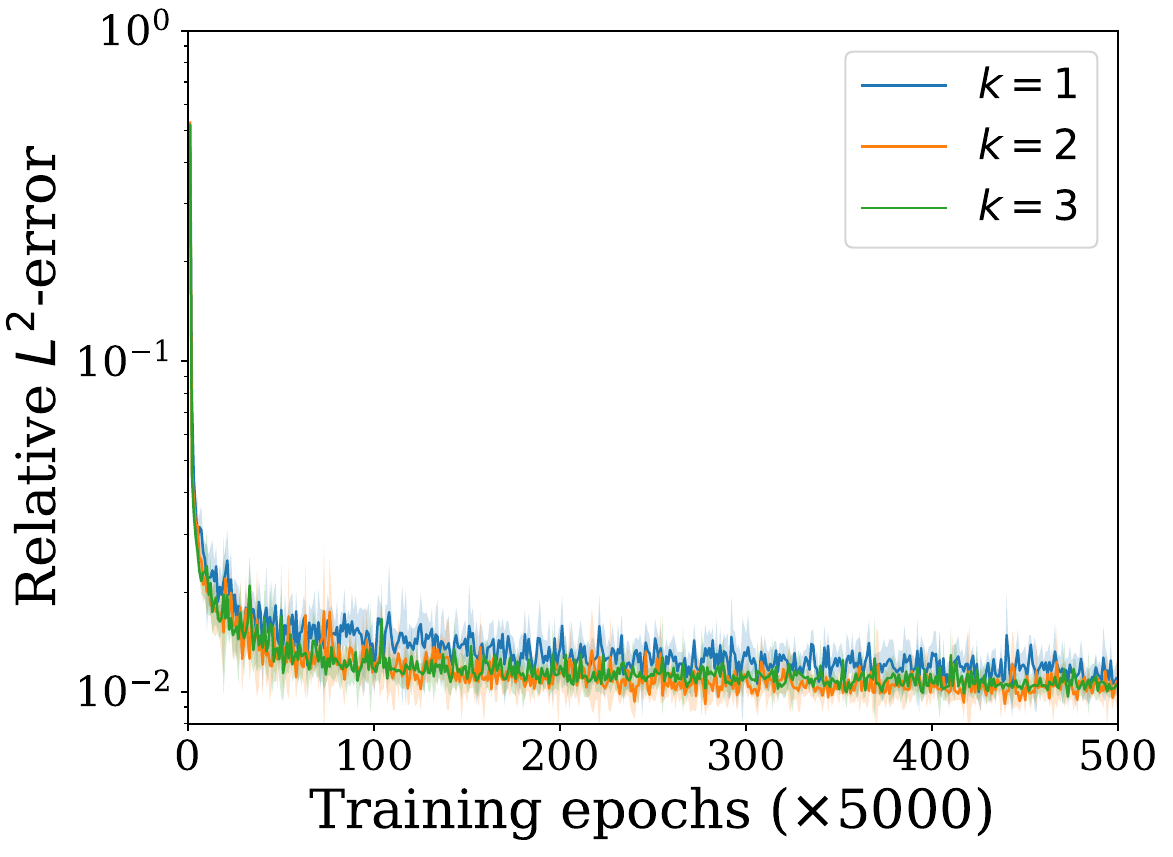}
        \caption{}
        \label{fig:err_convergence_3d_N:a}
    \end{subfigure}
    \hfill
    \begin{subfigure}[t]{0.49\textwidth}
        \centering
        \includegraphics[width=\textwidth]{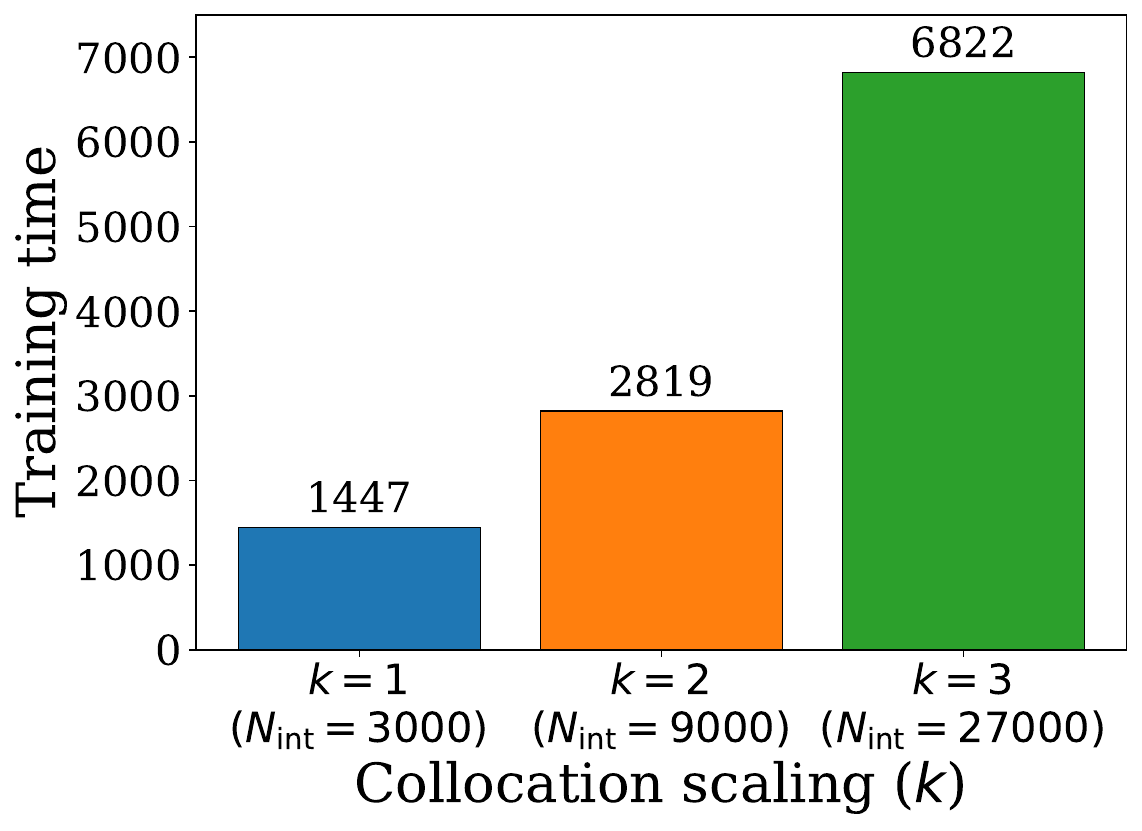}
        \caption{}
        \label{fig:err_convergence_3d_N:b}
    \end{subfigure}
    %% 5D
    \begin{subfigure}[t]{0.49\textwidth}
        \centering
        \includegraphics[width=\textwidth]{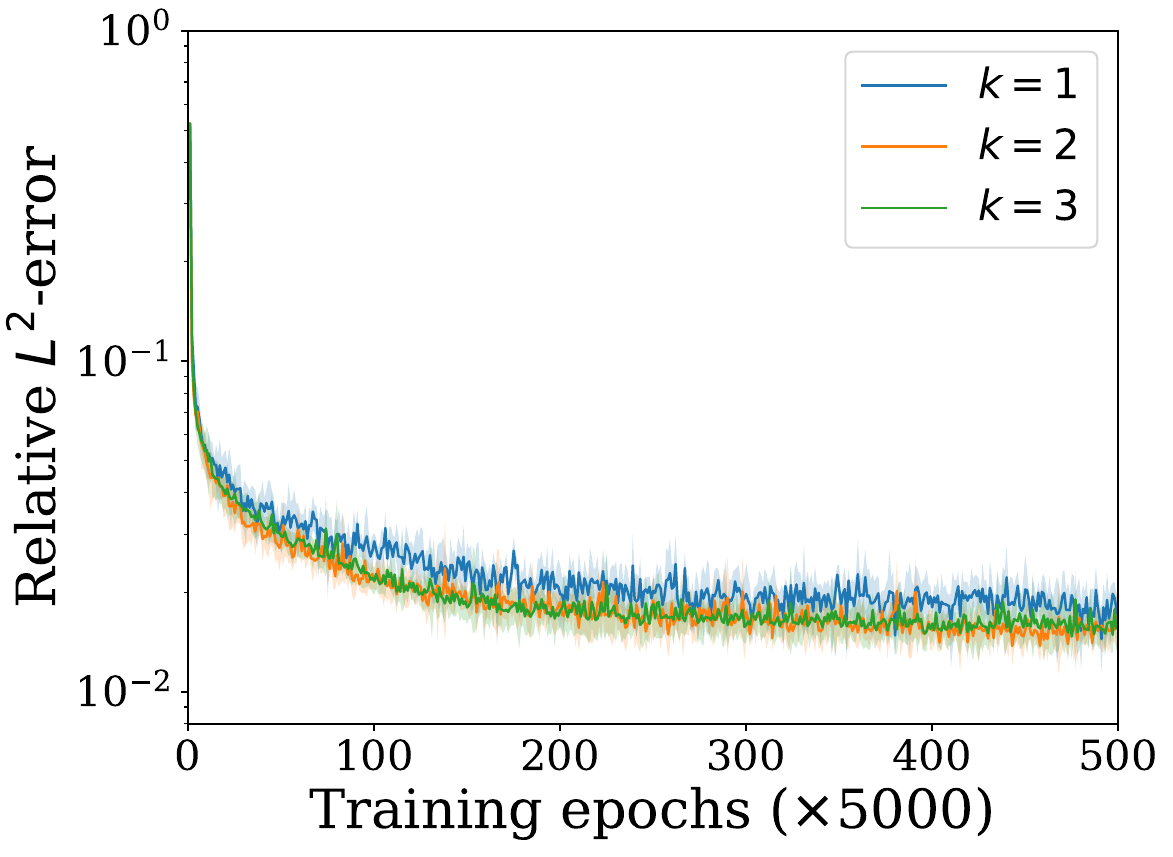}
        \caption{}
        \label{fig:err_convergence_3d_N:c}
    \end{subfigure}
    \hfill
    \begin{subfigure}[t]{0.49\textwidth}
        \centering
        \includegraphics[width=\textwidth]{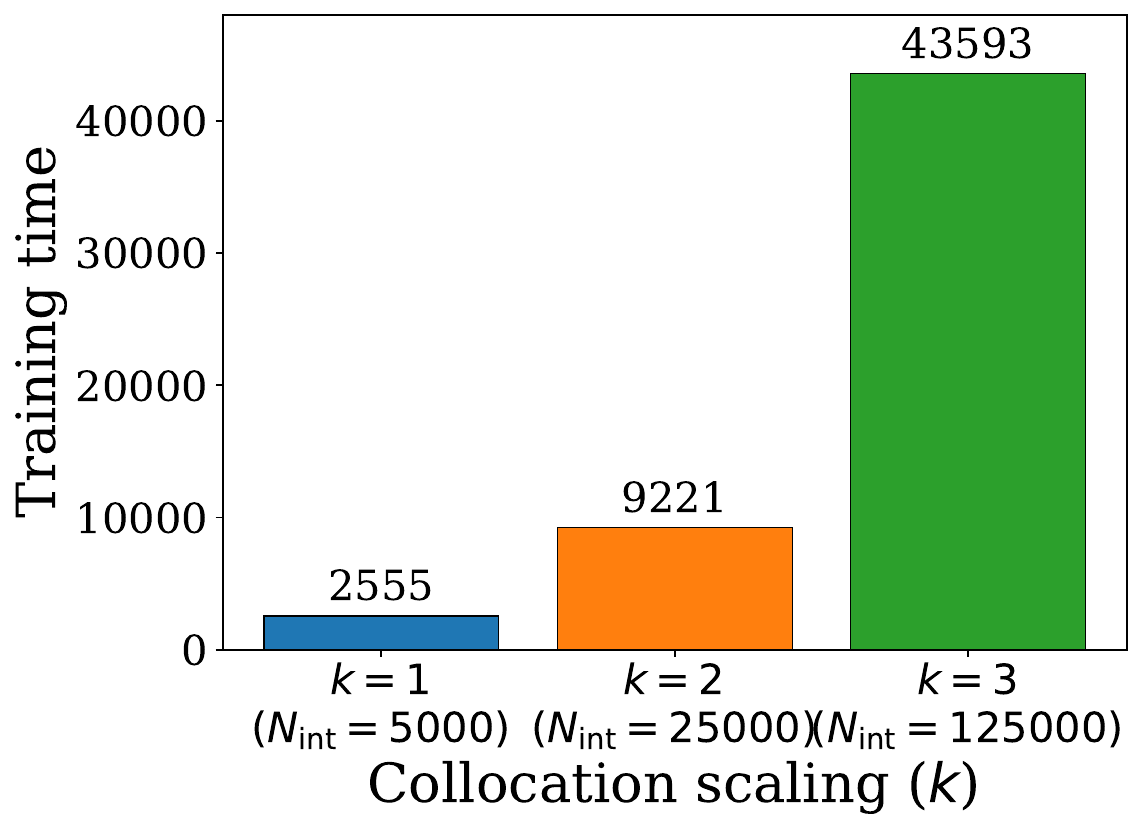}
        \caption{}
        \label{fig:err_convergence_3d_N:d}
    \end{subfigure}
    \caption{Effect of interior-collocation scaling on training dynamics.
    Top row: $N=1$ (3D). Bottom row: $N=2$ (5D).
    (a,c) Relative $L^2$-error versus training progress for $k\in\{1,2,3\}$.
    (b,d) Corresponding wall-clock training time for each $k$, where the number of interior points is
    $N_{\mathrm{int}} = 1000 (2N+1)^k$.}
    \label{fig:colloc_scaling}
\end{figure}

Finally, we examine the sensitivity of the policy-iterative PINN to the number of interior collocation points. Figure~\ref{fig:colloc_scaling} reports results for the publisher-subscriber problems with $N\in\{1,2\}$ (3D and 5D), where the number of interior points is scaled according to $N_{\mathrm{int}} = 1000 (2N+1)^k$ with $k\in\{1,2,3\}$. Increasing the collocation density improves accuracy, as expected, but the marginal gains beyond $k=2$ are limited relative to the additional computational cost. In particular, the configurations $k=2$ and $k=3$ exhibit similar convergence rates and comparable error floors.

Based on these observations, we fix $k=2$ in all subsequent experiments,
corresponding to $N_{\mathrm{int}} = 1000 (2N+1)^2$. This choice provides a favorable balance between accuracy and efficiency,
and suggests that the policy-iterative PINN is relatively insensitive
to further increases in collocation density once moderate domain coverage is achieved.

\subsection{Additional experiments: small diffusion and DeepONet comparison}\label{sec:additional_experiments}

We close the numerical study with two additional tests that complement the preceding experiments. The first examines the behavior of the policy-iterative PINN as the diffusion amplitude decreases and the ellipticity constant in the analysis becomes small. The second compares the method with a DeepONet operator-learning baseline on a family of moving-obstacle problems parameterized by the target location.

\paragraph{Small-diffusion regime.}
The convergence analysis in Section~\ref{sec:method} assumes a uniformly elliptic diffusion matrix, namely $\sigma\sigma^\top\succeq \lambda I$. For $\sigma=\gamma I_3$, the ellipticity constant scales as $\lambda=\gamma^2$. Thus, as $\gamma$ decreases, the publisher-subscriber problem approaches a weakly diffusive regime in which the second-order regularization is reduced. To examine this effect numerically, we consider the three-dimensional publisher-subscriber problem with
\[
\gamma\in\{0.10,0.07,0.05,0.03,0.01\}.
\]
The policy-iterative PINN and the direct PINN are trained with the same total optimization budget for each value of $\gamma$.

In this regime, the FDM reference also requires care. 
The central-difference computation used in the preceding experiments is performed on a grid with approximately $601^3$ spatial points, and the stability/monotonicity margin of this discretization deteriorates as the diffusion amplitude decreases. 
Maintaining the same effective resolution relative to the diffusion scale would require substantial mesh refinement in three dimensions. 
Indeed, since parabolic time stepping gives $\Delta t=O(h^2)$, halving the mesh size increases memory by approximately $2^3$ and the number of time steps by approximately $2^2$, resulting in an overall cost increase of roughly $2^5=32$. 
For this reason, we report FDM-relative errors only in regimes where the FDM reference remains sufficiently reliable, and use the smaller-diffusion profiles as qualitative diagnostics.

\begin{figure}[ht!]
    \centering
    \includegraphics[width=0.92\textwidth]{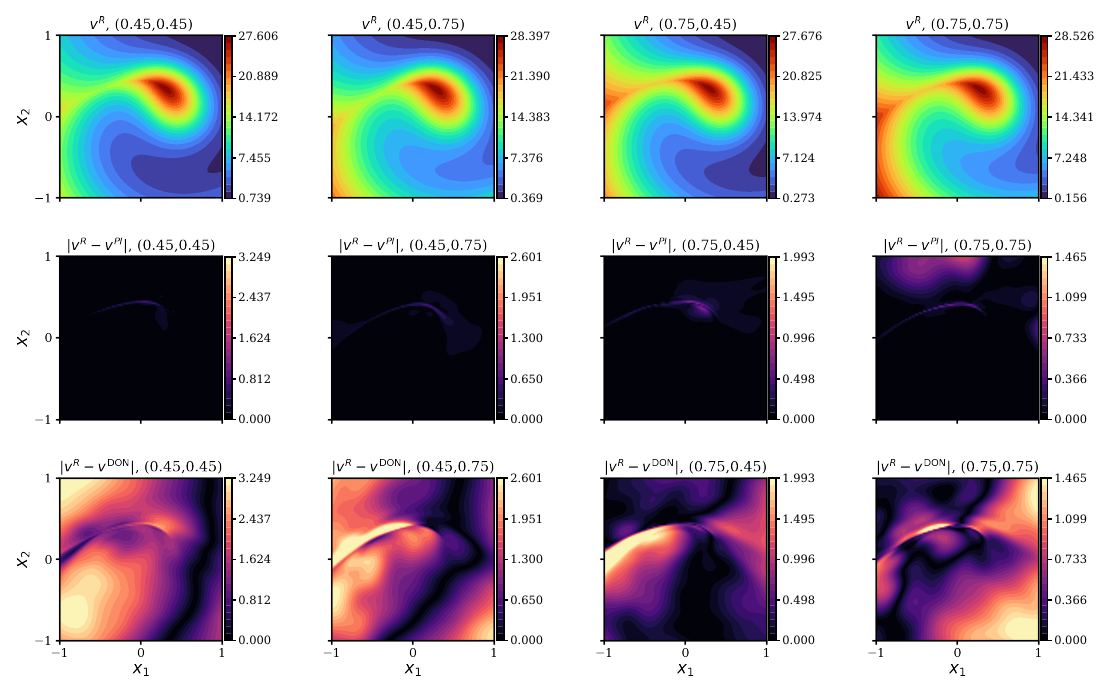}
    \caption{
    Small-diffusion test for the three-dimensional publisher-subscriber problem.
    Columns correspond to diffusion amplitudes $\gamma$ in $\sigma=\gamma I_3$.
    Rows show the FDM profile, the policy-iterative PINN solution, and the direct PINN solution at $t=0$ on the diagonal slice $x_1=x_2$.
    Within each column, the color scale is shared across the three rows.
    For $\gamma\le 0.05$, the FDM profiles are shown only as numerical profiles computed on the current grid, not as reliable quantitative references.
    }
    \label{fig:small_diffusion}
\end{figure}

\begin{table}[ht!]
    \caption{
    FDM status and estimated refinement cost for the small-diffusion publisher-subscriber experiment.
    The estimates indicate the cost of maintaining comparable effective resolution relative to the diffusion scale in three spatial dimensions.
    }
    \label{tab:small_diffusion_fdm_status}
    \begin{center}
        \resizebox{\textwidth}{!}{%
        \begin{tabular}{ccccc}
            \hline\hline
            $\sigma=\gamma I_3$ & $\gamma^2$ & FDM status with $h=1/600$ & estimated mesh spacing & estimated relative cost \\
            \noalign{\hrule height 1.5pt}
            $0.10I_3$ & $1.0\times 10^{-2}$ & reliable & current & $1$ \\
            $0.07I_3$ & $4.9\times 10^{-3}$ & borderline & current & $1$ \\
            $0.05I_3$ & $2.5\times 10^{-3}$ & not guaranteed & $\approx 1/1200$ & $\approx 32$ \\
            $0.03I_3$ & $9.0\times 10^{-4}$ & not guaranteed & $\approx 1/3300$ & $\approx 5.0\times 10^3$ \\
            $0.01I_3$ & $1.0\times 10^{-4}$ & not guaranteed & $\approx 1/30000$ & prohibitive \\
            \noalign{\hrule height 1.5pt}
        \end{tabular}
        }
    \end{center}
    \vskip-.2truecm
\end{table}

\begin{table}[ht!]
    \caption{Quantitative error comparison for the small-diffusion publisher-subscriber experiment. Relative $L^2$-errors at $t=0$ are reported only for the diffusion levels where the current FDM profile is used as a reliable or borderline reference.}
    \label{tab:small_diffusion_error_comparison}
    \begin{center}
        \begin{tabular}{cccc}
            \hline\hline
            $\sigma=\gamma I_3$ & FDM status & PI-PINN & Direct PINN \\
            \noalign{\hrule height 1.5pt}
            $0.10I_3$ & reliable & $8.947{\times}10^{-3}$ & $1.200{\times}10^{-2}$ \\
            $0.07I_3$ & borderline & $6.450{\times}10^{-3}$ & $1.176{\times}10^{-2}$ \\
            \noalign{\hrule height 1.5pt}
        \end{tabular}
    \end{center}
    \vskip-.2truecm
\end{table}

For the two diffusion levels for which quantitative FDM comparison is reported, the policy-iterative PINN attains smaller final relative $L^2$-errors than the direct PINN under the same optimization budget; see Table~\ref{tab:small_diffusion_error_comparison}. Figure~\ref{fig:small_diffusion} gives the corresponding solution profiles across all five diffusion amplitudes. For $\gamma\le 0.05$, however, we do not use the current central-difference profile for quantitative FDM benchmarking because its monotonicity-based reliability is not guaranteed, as summarized in Table~\ref{tab:small_diffusion_fdm_status}. These profiles therefore serve as a qualitative diagnostic of the neural solvers as the second-order regularization weakens.

\paragraph{Comparison with a DeepONet baseline.}
We further compare the proposed policy-iterative PINN with a DeepONet baseline on the two-dimensional moving-obstacle problem. 
Since DeepONet is designed to learn a solution operator over a family of problem instances, we parameterize the family by the target location $x_{\mathrm{goal}}$. 
The DeepONet baseline is trained once on the $3\times3$ target-location grid
\[
\{0.3,0.6,0.9\}\times\{0.3,0.6,0.9\},
\]
and is evaluated on the held-out target locations
\[
x_{\mathrm{goal}}\in\{(0.45,0.45),(0.45,0.75),(0.75,0.45),(0.75,0.75)\}.
\]
In contrast, the policy-iterative PINN is trained independently for each target location. 
Both methods are evaluated against the same FDM reference solutions using the relative $L^2$-error at $t=0$.

\begin{figure}[ht!]
    \centering
    \includegraphics[width=0.92\textwidth]{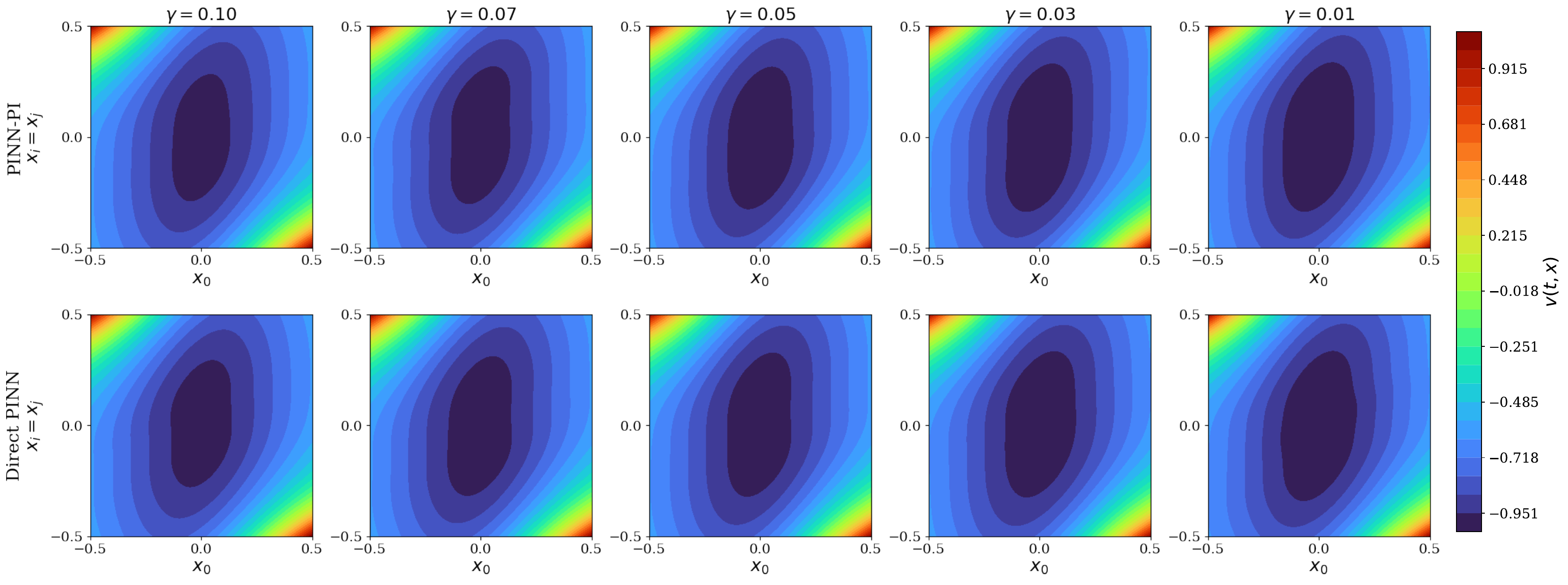}
    \caption{
    Target-location comparison for the two-dimensional moving-obstacle problem at $t=0$.
    Columns correspond to held-out target locations.
    The first row shows the FDM reference value, and the second and third rows show the absolute errors of the policy-iterative PINN and DeepONet (DON), respectively.
    }
    \label{fig:deeponet_target}
\end{figure}

For the DeepONet baseline, the branch input consists of terminal-cost values evaluated at $81$ fixed sensor points, while the trunk input is $(t,x)$.
We use $2{,}000$ spatio-temporal collocation points per training target, resulting in $18{,}000$ collocation points per outer iteration. 
For the policy-iterative PINN, we report the mean, standard deviation, and best final relative $L^2$-error over five independent random seeds. 
For DeepONet, we report the mean, standard deviation, and best relative $L^2$-error over five independently trained operators evaluated on the held-out targets.
The quantitative errors are summarized in Table~\ref{tab:deeponet_comparison}.
Figure~\ref{fig:deeponet_target} visualizes, for each target location, the best-performing PI-PINN and DeepONet runs reported in Table~\ref{tab:deeponet_comparison}.

\begin{table}[ht!]
    \caption{
    Target-location comparison for the two-dimensional moving-obstacle problem.
    Errors are relative $L^2$-errors at $t=0$ against the FDM reference.
    PI-PINN rows report the mean, standard deviation, and best value over five independent random seeds.
    DeepONet rows report the same statistics over five independently trained operators, each trained on the $3\times3$ target-location grid and evaluated on the held-out targets.
    }
    \label{tab:deeponet_comparison}
    \begin{center}
        \resizebox{\textwidth}{!}{%
        \begin{tabular}{cccccc}
            \hline\hline
            Method & Metric 
            & $(0.45,0.45)$ 
            & $(0.45,0.75)$ 
            & $(0.75,0.45)$ 
            & $(0.75,0.75)$ \\
            \noalign{\hrule height 1.5pt}
            \multirow{3}{*}{PI-PINN}
            & Mean 
            & $1.525{\times}10^{-2}$ 
            & $8.777{\times}10^{-3}$ 
            & $9.226{\times}10^{-3}$ 
            & $2.440{\times}10^{-2}$ \\
            & Std 
            & $1.121{\times}10^{-2}$ 
            & $2.884{\times}10^{-3}$ 
            & $2.999{\times}10^{-3}$ 
            & $6.966{\times}10^{-3}$ \\
            & Best 
            & $6.135{\times}10^{-3}$ 
            & $6.454{\times}10^{-3}$ 
            & $6.695{\times}10^{-3}$ 
            & $1.270{\times}10^{-2}$ \\
            \noalign{\hrule height 1.5pt}
            \multirow{3}{*}{DeepONet}
            & Mean
            & $5.648{\times}10^{-1}$
            & $2.428{\times}10^{-1}$
            & $3.437{\times}10^{-1}$
            & $1.940{\times}10^{-1}$ \\
            & Std
            & $2.417{\times}10^{-1}$
            & $7.810{\times}10^{-2}$
            & $1.255{\times}10^{-1}$
            & $4.333{\times}10^{-2}$ \\
            & Best
            & $9.263{\times}10^{-2}$
            & $1.320{\times}10^{-1}$
            & $1.078{\times}10^{-1}$
            & $1.342{\times}10^{-1}$ \\
            \noalign{\hrule height 1.5pt}
        \end{tabular}
        }
    \end{center}
    \vskip-.2truecm
\end{table}

Figure~\ref{fig:deeponet_target} and Table~\ref{tab:deeponet_comparison} illustrate the different computational tradeoffs of the two approaches.
DeepONet amortizes training over a family of target locations and evaluates held-out targets through a single learned operator.
By contrast, the policy-iterative PINN solves each target instance separately through residual minimization under policy updates.
In the present experiment, the per-instance PI-PINN achieves smaller FDM-relative errors at $t=0$, while the DeepONet result reflects the accuracy of an amortized operator approximation on held-out target locations.
This comparison should therefore be interpreted as a tradeoff between amortized operator learning across related problem instances and dedicated residual-based optimization for each fixed game configuration.

\section{Conclusion}\label{sec:conclusion}
In this work, we have proposed a novel mesh-free framework for solving nonconvex Hamilton--Jacobi--Isaacs (HJI) equations by combining classical policy iteration with physics-informed neural networks (PINNs). By leveraging the differentiability and expressive capacity of deep neural networks, the proposed method enables efficient approximation of viscosity solutions to high-dimensional HJI equations, even in the presence of nonconvexities in the Hamiltonian. We have provided a rigorous convergence analysis under a uniform ellipticity condition and demonstrated the numerical effectiveness of the method across several benchmark differential games in two to eleven dimensions.

The experiments confirm that our approach not only achieves competitive accuracy but also exhibits improved convergence behavior as the problem dimension increases, compared to standard PINN and direct collocation methods. Moreover, empirical results show that the iterative policy-improvement scheme yields lower residuals than direct one-shot training and achieves faster convergence and lower error floors in isotropic settings, and remains stable in highly nonconvex regimes. This is attributed to the fact that each policy-evaluation step solves a linearized PDE with fixed control inputs, resulting in smoother and more stable optimization dynamics. As a result, the learned value functions are not only more accurate but also demonstrate improved accuracy and more stable convergence, even in highly nonconvex settings.

Furthermore, unlike direct one-shot approaches, our iterative framework enables systematic quantification of approximation errors. In particular, the $L^2$-error between the learned and true value functions can be rigorously bounded even in nonconvex settings, owing to the linearized structure of each policy evaluation step. Nevertheless, the method inherits limitations intrinsic to PINN-based approaches, such as sensitivity to network initialization and challenges in gradient stability over long training horizons.

An important limitation of the present work is the assumption of non-degenerate diffusion. In the absence of noise, the problem becomes a first-order Hamilton--Jacobi--Isaacs equation arising in deterministic differential games. This setting lacks the regularizing properties of second-order terms and poses distinct theoretical and computational challenges. We leave the development of PINN-based policy iteration schemes for such first-order HJI problems to future work. Nonconvex HJ equations with degeneracy have seen significant progress, notably through policy iteration~\cite{guo2025policy} and the nonlinear adjoint approach~\cite{mitake2016dynamical}. We expect that combining these ideas with our method could lead to effective extensions to degenerate or first-order problems.

Future work will aim to address these issues by incorporating adaptive sampling strategies, exploring PINN variants based on operator learning (e.g., DeepONet, Fourier Neural Operator), and extending the framework to stochastic differential games and time-varying dynamics. The integration of our PINN-based policy iteration with model-based reinforcement learning paradigms also presents an exciting direction for further research.

\section*{Acknowledgement}
Yeoneung Kim is supported by the National Research Foundation of Korea (NRF) grant funded by the Korea government(MSIT) (RS-2023-00219980, RS-2026-25490291). Hee Jun Yang and Minjung Gim are supported by National Institute for Mathematical Sciences (NIMS) grant funded by the Korea government (MSIT) (No. B26810000). The authors would like to thank Professor Hung Vinh Tran (University of Wisconsin-Madison) for his insightful suggestions and valuable guidance in developing and refining the ideas of this work.

\section*{Conflict of interest}
The authors declare that they have no known competing financial interests or personal relationships that could have appeared to influence the work reported in this paper.

\section*{Data availability}
The datasets generated during and/or analyzed during the current study are available in the GitHub repository, \url{https://github.com/HeeJunYang09/hji-pi-pinn-repro}.

\section*{CRediT authorship contribution statement}
	\textbf{Hee Jun Yang:} Software, Investigation, Validation, Visualization, Writing - original draft, Writing - review \& editing
	\textbf{Minjung Gim:} Writing - original draft, Writing - review \& editing, Supervision
	\textbf{Yeoneung Kim:} Conceptualization, Methodology, Formal analysis, Writing - original draft, Writing - review \& editing, Supervision

\bibliographystyle{plain}
\bibliography{ref.bib}

\appendix

\section{Parabolic $L^2$-estimate}\label{app:parabolic}

\begin{prop}[Parabolic $L^{2}$-estimate]\label{prop:L2}
Let \(d\in\mathbb N\) and \(T>0\), and suppose
Assumption~\ref{assump:main} holds. Let
\[
b\in L^\infty((0,T)\times\mathbb R^d;\mathbb R^d),
\qquad
P\in L^2(0,T;L^2(\mathbb R^d)),
\qquad
Q\in L^2(\mathbb R^d).
\]
If
\[
v\in C([0,T];L^2(\mathbb R^d))
\cap L^2(0,T;H^1(\mathbb R^d)),
\qquad
\partial_t v\in L^2(0,T;H^{-1}(\mathbb R^d)),
\]
is a variational solution of
\begin{equation}\label{eq:PE_L2_SIAM}
\begin{cases}
\displaystyle
\partial_t v
+\frac12\operatorname{Tr}
\bigl(\sigma\sigma^\top D_{xx}^2v\bigr)
+b(t,x)\cdot\nabla_xv
=P(t,x),
&
\textnormal{in }[0,T)\times\mathbb R^d,
\\[0.4em]
v(T,\cdot)=Q,
&
\textnormal{on }\mathbb R^d,
\end{cases}
\end{equation}
then there exists
\[
C_T
=
C\!\left(
T,\lambda,\|b\|_\infty,
\|\operatorname{div}(\sigma\sigma^\top)\|_\infty
\right)>0
\]
such that
\[
\sup_{t\in[0,T]}\|v(t,\cdot)\|_2^2
+
\lambda\int_0^T
\|\nabla_xv(t,\cdot)\|_2^2\,\de t
\le
C_T
\left(
\|Q\|_2^2
+
\int_0^T\|P(t,\cdot)\|_2^2\,\de t
\right).
\]
\end{prop}

\begin{proof}
This is the standard parabolic energy estimate; see, e.g., \cite[Section~7.1]{evans2022partial} and \cite[Chapter~III]{ladyzhenskaya1968linear} for the classical $L^2$-theory. Testing
\eqref{eq:PE_L2_SIAM} with \(2v\), using
\[
\operatorname{Tr}(AD^2v)
=
\operatorname{div}(A\nabla v)
-
(\operatorname{div}A)\cdot\nabla v,
\qquad A=\sigma\sigma^\top,
\]
and applying uniform ellipticity, Young's inequality, and Gronwall's
lemma yield the result.
\end{proof}

\section{Nonconvexity of the Hamiltonians in the numerical examples}
\label{app:nonconvex_example}
In this appendix, we verify that the Hamiltonians arising in the main numerical examples of this paper are genuinely nonconvex with respect to the gradient variable.

\subsection{Two-dimensional optimal path planning with a moving obstacle}
We recall the Hamiltonian associated with the two-dimensional path planning problem introduced in Section~\ref{sec:pathplanning}.
It is given by
\[
H(x,p)
=
\sup_{|b|\le \delta}
\inf_{|a|\le 1}
\left(
\lambda_1 |a|^2
+ \lambda_2 \phi(x)
+ p\cdot(a+b)
\right).
\]

Evaluating the minimax problem yields the explicit form
\[
H(x,p)
=
\begin{cases}
-\dfrac{1}{4\lambda_1}|p|^2 + \lambda_2 \phi(x) + \delta |p|,
& |p|\le 2\lambda_1,\\[0.5em]
-|p| + \lambda_1 + \lambda_2 \phi(x) + \delta |p|,
& |p|>2\lambda_1.
\end{cases}
\]

While each individual term is convex or concave in isolation,
the combination of the quadratic term $-|p|^2$ and the linear absolute-value term $\delta |p|$ induces a Hamiltonian that is neither convex nor concave in $p$. In particular, the change of regime at $|p|=2\lambda_1$
creates a non-smooth and nonconvex geometry in the gradient variable.

\subsection{Publisher-subscriber problem with $N=1$}

The Hamiltonian associated with the publisher-subscriber game takes the form
$$
H(x,p)
=
p^\top (Ax+\psi(x))
-\|B^\top p\|_1
+\|C^\top p\|_1,
\qquad \|z\|_1:=\sum_k |z_k|,
$$
which arises from bounded control and disturbance sets.
The linear term $p^\top (Ax+\psi(x))$ does not affect convexity in $p$, so the (non)convexity is determined by the Isaacs part
\begin{equation*}
\mathcal{I}(p):=-\|B^\top p\|_1+\|C^\top p\|_1.
\end{equation*}
When the disturbance channel is axis-aligned, the two $\ell_1$ terms in $\mathcal{I}(p)=-\|B^\top p\|_1+\|C^\top p\|_1$ are aligned in the same coordinate directions, so they effectively combine into a single $\ell_1$-norm term (up to a constant factor), which is convex (or, with a negative sign, concave).
With a rotated disturbance map, the second $\ell_1$ norm is taken after a rotation, so its kinks occur along different directions than those of the first term; this misalignment prevents cancellation and typically makes $\mathcal{I}(p)$ neither convex nor concave.

Consider the three-dimensional instance with one subscriber pair, so that the state dimension is $3$ and $p=[p_1,p_2,p_3]^\top\in\mathbb{R}^3$.
In this case,
\begin{equation*}
B=\begin{bmatrix}\mathbf 0_{1\times 2}\\ b I_2\end{bmatrix},
\qquad
C=\begin{bmatrix}\mathbf 0_{1\times 2}\\ c R(\phi)\end{bmatrix},
\qquad
R(\phi)=
\begin{bmatrix}
\cos\phi & -\sin\phi\\
\sin\phi & \cos\phi
\end{bmatrix},
\end{equation*}
and thereby,

\begin{equation*}
B=\begin{bmatrix}
0 & 0\\
b & 0\\
0 & b
\end{bmatrix},
\qquad
C=\begin{bmatrix}
0 & 0\\
c\cos\phi & -c\sin\phi\\
c\sin\phi & \ \ c\cos\phi
\end{bmatrix}.
\end{equation*}
Hence, 
\[
B^\top p=\begin{bmatrix} bp_2\\ b p_3\end{bmatrix},
\qquad
C^\top p=c R(\phi)^\top\begin{bmatrix} p_2\\  p_3 \end{bmatrix},
\]
and the Isaacs term reduces to a function of the subscriber-gradient components:
\begin{equation}
\mathcal{I}(p)
=
-\|B^\top p\|_1
+\|C^\top p\|_1.
\label{eq:Isaacs_N1_app}
\end{equation}
% \begin{equation}
% \mathcal{I}(p)
% =
% -\|b(p_2,p_3)\|_1
% +c\|R(\phi)^\top (p_2,p_3)\|_1.
% \label{eq:Isaacs_N1_app}
% \end{equation}

Each norm in \eqref{eq:Isaacs_N1_app} is convex, but their difference is not convex in general. Moreover, since the $\ell_1$ norm is not rotation-invariant, the two terms have kinks along different directions in the $(p_2,p_3)$-plane whenever $\phi\notin\{0,\pi\}$.
As a result, $\mathcal{I}(p)$ is generally neither convex nor concave in $p$ (with the axis-aligned cases $\phi=0$ or $\pi$ being the main exceptions, where the kink sets align and $\mathcal{I}(p)$ simplifies to a single $\ell_1$ term up to a scalar).

% \begin{figure}[ht!]
%   \centering
%   \includegraphics[width=\textwidth]{hamiltonian.pdf}
%   \caption{
%   Visualization of $H(x,p)$ in the minimal case $N=1$ (one publisher, two subscribers) with rotation angle $\phi=\pi/4$.
%   We fix $x=[0,0,0]^\top$ and plot $H(x,p)$ over three 2D slices: $[p_1,p_2]^\top$ with $p_3=0$ (left), $[p_1,p_3]^\top$ with $p_2=0$ (middle), and $[p_2,p_3]^\top$ with $p_1=0$ (right), each over $[-2,2]^2$.
%   The rotated disturbance map induces kink directions that differ from those of the control channel, producing a piecewise-linear geometry that is generally neither convex nor concave in $p$.
%   }
%   \label{fig:hamiltonian_nonconvex_app}
% \end{figure}

\begin{figure}[ht!]
  \centering
  \begin{subfigure}[b]{0.50\textwidth}
    \centering
    \includegraphics[width=\textwidth]{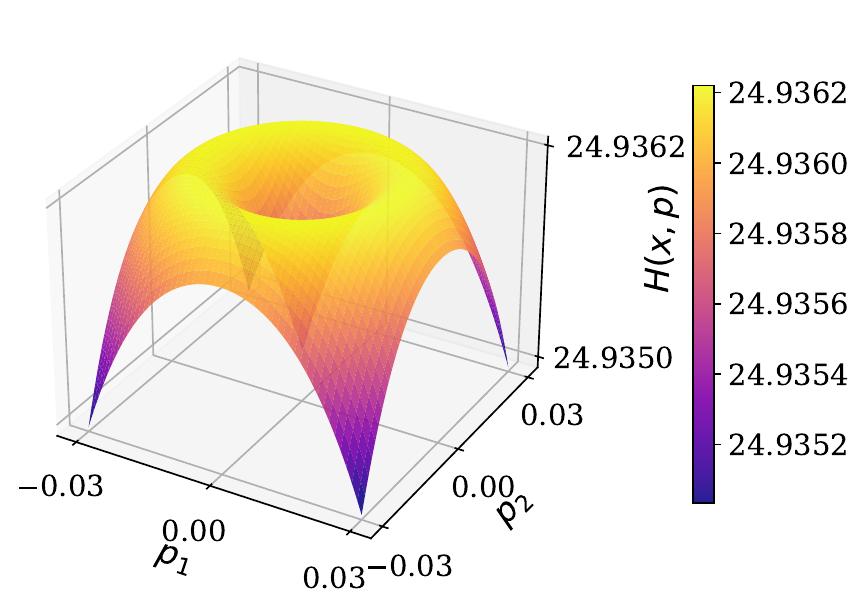}
    \caption{}
    \label{fig:ham_path}
  \end{subfigure}
  \hfill
  \begin{subfigure}[b]{0.46\textwidth}
    \centering
    \includegraphics[width=\textwidth]{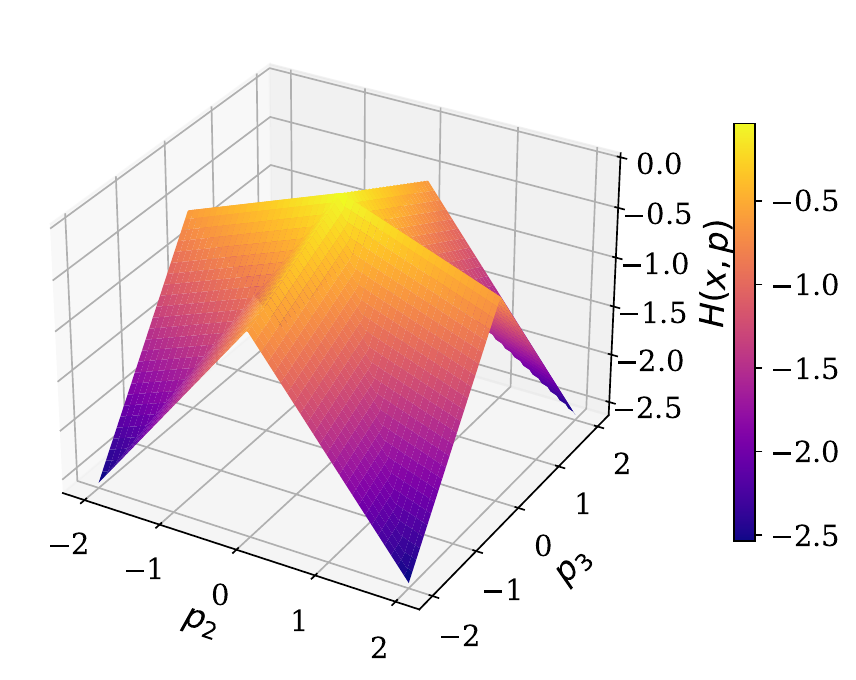}
    \caption{}
    \label{fig:ham_pubsub}
  \end{subfigure}
  
  \caption{
  Comparison of Hamiltonian visualizations. (a) $H(x,p)$ for the path planning problem. (b) $H(x,p)$ for the publisher-subscriber problem over the $[p_2, p_3]^\top$ slice.
  }
  \label{fig:hamiltonian_nonconvex_app}
\end{figure}

\section{Proof of Value Function Decomposition}\label{proof:value_decomp}

We show that the value function $v(t,x)$ associated with the $(2N+1)$-dimensional HJI equation admits an exact decomposition across publisher-subscriber pair blocks, provided that (i) the terminal cost is separable across blocks, (ii) the drift is unidirectionally coupled (publisher $\to$ subscribers) in the sense specified below, and (iii) the diffusion is block-diagonal consistent with the same partition.

Let \(x=[x_0,x_1,\ldots,x_{2N}]^\top\in\mathbb{R}^{2N+1}\), where \(x_0\) is the publisher state and
\([x_{2i-1},x_{2i}]^\top\in\mathbb{R}^2\) is the \(i\)th subscriber pair for \(i=1,\ldots,N\).
We assume that the terminal cost is separable across the \(N\) publisher-subscriber pair blocks:
\[
g(x)=\sum_{i=1}^{N} g_i(P_i x),
\qquad
g_i(z):=\tfrac12\left(z_1^2+z_2^2+z_3^2-r^2\right),\quad z\in\mathbb{R}^3,
\]
where \(P_i:\mathbb{R}^{2N+1}\to\mathbb{R}^3\) denotes the projection onto the coordinates
\((x_0,x_{2i-1},x_{2i})\), i.e.,
\[
P_i x:=\left[x_0,  x_{2i-1},  x_{2i}\right]^\top.
\]
Equivalently, \(g_i(P_i x)=\tfrac12\left(x_0^2+x_{2i-1}^2+x_{2i}^2-r^2\right)\) for each \(i\). Let $v$ be the (unique) viscosity solution of the full HJI equation

\[
\begin{cases}
	\partial_t v(t,x) + H(t,x,\nabla_x v(t,x)) = -\frac{1}{2} \operatorname{Tr}(\sigma\sigma^\top D_{xx}^2 v(t,x)) \quad&\text{in}\quad [0,T) \times \mathbb{R}^{2N+1},\\
	v(T,x)=g(x) \quad&\text{on}\quad \mathbb{R}^{2N+1}.
\end{cases}
\]

We define the candidate decomposed value function
$$
\hat v(t,x):=\sum_{i=1}^N v_i\left(t,x_0,x_{2i-1},x_{2i}\right),
$$
where each $v_i$ is a viscosity solution of a three-dimensional reduced HJI equation in the variables $[x_0,x_{2i-1},x_{2i}]^\top$, treating all other state coordinates as absent:
\[
\begin{cases}
\partial_t v_i + H_i\left(t, x_0,x_{2i-1},x_{2i},\nabla_x v_i\right)
= -\dfrac12 \operatorname{Tr}\left(\sigma_i\sigma_i^\top D^2_{xx} v_i\right) \quad&\text{in}\quad [0,T) \times \mathbb{R}^{3},\\
	v_i(T,P_i x)=g_i(P_i x) \quad&\text{on}\quad \mathbb{R}^{3},
\end{cases}
\]
with $\nabla_x v_i=[\partial_{x_0}v_i,\partial_{x_{2i-1}}v_i,\partial_{x_{2i}}v_i]^\top$.

Here $H_i$ denotes the reduced Hamiltonian obtained by restricting the full dynamics and admissible inputs to the $i$th block. Concretely, the unidirectional coupling assumption is that the publisher coordinate $x_0$ evolves independently of subscriber coordinates, and each subscriber pair $[x_{2i-1},x_{2i}]^\top$ depends on $[x_0,x_{2i-1},x_{2i}]^\top$ but not on other subscriber pairs.

We assume the diffusion matrix is block-diagonal with respect to the partition $x_0$ and $[x_{2i-1},x_{2i}]^\top$:
$$
\sigma=\mathrm{diag}\left(\sigma_0,\Sigma_1,\ldots,\Sigma_N\right)\in \mathbb{R}^{(2N+1)\times (2N+1)},
$$
where $\sigma_0\in\mathbb{R}$ and $\Sigma_i\in\mathbb{R}^{2\times2}$. Let $\sigma_i\in\mathbb{R}^{3\times 3}$ denote the induced diffusion on the $i$th $[x_0,x_{2i-1},x_{2i}]^\top$ block, i.e.,
$$
\sigma_i:=\mathrm{diag}(\sigma_0,\Sigma_i).
$$
The following calculation is understood in the viscosity sense; it can
be justified rigorously by the standard parabolic theorem of sums
(cf.~\cite[Theorem~8.3]{crandall1992user}).

Since each $v_i$ depends on $(x_0,x_{2i-1},x_{2i})$, the mixed second derivatives across different subscriber blocks vanish. Hence the second-order term decomposes as
$$
\operatorname{Tr}\left(\sigma\sigma^\top D^2_{xx}\hat v\right)
=\sum_{i=1}^N \operatorname{Tr}\left(\sigma_i\sigma_i^\top D_{xx}^2 v_i\right),
$$
where $D^2_{xx} v_i$ is the $3\times 3$ Hessian with respect to $[x_0,x_{2i-1},x_{2i}]^\top$.

Because $\hat v$ is a sum of block-wise functions, its derivatives decompose as
$$
\partial_t \hat v=\sum_{i=1}^N \partial_t v_i,\qquad
\nabla_x \hat v=\sum_{i=1}^N E_i \nabla_x v_i,
$$
where $E_i:\mathbb{R}^3\to\mathbb{R}^{2N+1}$ embeds a 3-vector into the full space by placing entries in coordinates $(0,2i-1,2i)$ and zeros elsewhere. 

Under the unidirectional coupling, the full Hamiltonian evaluated at $\nabla_x\hat v$ separates across blocks:
$$H\left(t,x,\nabla_x \hat v(t,x)\right)
=\sum_{i=1}^N H_i\left(t,x_0,x_{2i-1},x_{2i},\nabla_x v_i(t,x_0,x_{2i-1},x_{2i})\right).
$$

Intuitively, this holds because the drift and the control/disturbance channels act block-wise on subscriber pairs (with the shared publisher coordinate), so the dependence on $\nabla_x \hat v$ enters only through the corresponding block gradients. In particular, rotated disturbance maps of the form 
$$
C = \begin{bmatrix} \mathbf 0 \\ c \text{diag}(R(\phi_1),\ldots,R(\phi_N)) \end{bmatrix}.
$$
preserve this block-wise structure. Combining the decompositions of $\partial_t\hat v$, the Hamiltonian, and the second-order term, we conclude that $\hat v$ satisfies the same HJI equation as $v$.

By construction and separability of $g$, 
$$
\hat v(T,x)=\sum_{i=1}^N v_i(T,x_0,x_{2i-1},x_{2i})
=\sum_{i=1}^N g_i(x_0,x_{2i-1},x_{2i})
=g(x).
$$
Therefore, $\hat v$ is a viscosity solution to the full HJI terminal value problem with the same terminal data. Under standard assumptions ensuring uniqueness of viscosity solutions (e.g., Lipschitz continuity of the Hamiltonian in $p$ and uniform ellipticity of $\sigma\sigma^\top)$, it follows that
$$
v(t,x)=\hat v(t,x),\quad \text{for all }(t,x)\in[0,T]\times\mathbb{R}^{2N+1},
$$
and the exact decomposition holds globally.

\section{Details of Implementation Parameters}\label{sec:impl_details}
\begin{table}[H]
	\centering
	\caption{Summary of numerical settings for each experiment.}
	\label{tab:numerical_settings}
	\resizebox{\textwidth}{!}{%
		\begin{tabular}{lcc}
			\hline
			\textbf{Setting} & \textbf{Moving Obstacle (2D)} & \textbf{Publisher-Subscriber ($(2N+1)$D)} \\
			\hline
			\multicolumn{3}{l}{\textbf{Neural Network Settings}} \\
			Network architecture & $4$ hidden layers, 64 units & $4$ hidden layers, 64 units \\
			Training epochs per iteration ($E$) & $1{,}000$ & $5{,}000$ \\
			Policy iteration ($M$) & $1{,}000$ & $500$ \\
			Direct PINN epoch ($E\times M$) &         & $2{,}500{,}000$ \\
			Collocation points & $2{,}000$ (refreshed every 100 epochs) & $1{,}000 {(2N+1)}^2$ (refreshed every 500 epochs) \\
			Initial policy & Uniform over admissible set & Uniform over admissible set \\
			Extended spatial domain & $[-1, 1]^2$ & $[-1.5, 1.5]^{2N+1}$ \\
			Target domain & $[-1, 1]^2$ & $[-0.5, 0.5]^{2N+1}$ \\
			\hline
			\multicolumn{3}{l}{\textbf{Reference Solution Settings}} \\
			Extended spatial domain & $[-2, 2]^2$ & $[-1.5, 1.5]^3$ \\
			Target domain & $[-1, 1]^2$ & $[-0.5, 0.5]^3$ \\
			FDM spatial grid & $201 \times 201$ & $601 \times 601 \times 601$ \\
			FDM time steps & $201^2$ & $601^2$ \\
			Boundary condition & Homogeneous Neumann & Homogeneous Neumann \\
			\hline
		\end{tabular}
	}
    {\footnotesize $^{*}$Publisher-Subscriber: PINN solves the $(2N+1)$D problem directly; the reference is constructed from 3D projected subproblems via decomposition (Appendix~\ref{proof:value_decomp}).}
\end{table}

\end{document}